\documentclass[11pt, reqno, psamsfonts]{amsart}
\pdfoutput=1

\usepackage{amssymb}
\usepackage{amsthm}
\usepackage{amsmath}
\usepackage{latexsym}
\usepackage[T1]{fontenc}
\usepackage[utf8]{inputenc}
\usepackage[russian, french, english]{babel}

\usepackage{graphicx}
\usepackage{wrapfig}
\usepackage[justification=centering, labelfont=bf]{caption}
\usepackage{subcaption}
\usepackage{mathtools}
\usepackage[hidelinks]{hyperref}
\usepackage{tikz}
\usetikzlibrary{shapes,snakes}
\usetikzlibrary{arrows.meta}
\usepackage{amsbsy}
\usepackage{enumitem}
\usepackage{mathrsfs}
\usepackage{multicol}
\usepackage{bm}
\usepackage{thmtools}
\usepackage{marginnote}
\usepackage[T1]{fontenc}

\usepackage{Baskervaldx} 
\usepackage[baskervaldx,vvarbb]{newtxmath} 
\usepackage[bb=px,cal=boondoxo]{mathalfa}

\usepackage[spacing=true,kerning=true,babel=true,tracking=true]{microtype}

\usepackage[shortcuts]{extdash}

\usepackage[foot]{amsaddr} 

\usepackage{framed}

\usepackage[multiple]{footmisc}

\usepackage[left=1in,right=1in,top=1in,bottom=1in,bindingoffset=0cm]{geometry}

\usepackage[backend=biber, style=alphabetic, sorting=nyt, maxnames=100]{biblatex}
\title{Independent Sets in Algebraic Hypergraphs}
\date{}
\author{Anton Bernshteyn}
\address[Anton Bernshteyn]{\normalfont{}Department of Mathematical Sciences, Carnegie Mellon University, Pittsburgh, PA,	USA}
\email{abernsht@math.cmu.edu}

\author{Michelle Delcourt}
\address[Michelle Delcourt]{\normalfont{}Combinatorics and Optimization Department, University of Waterloo, Waterloo, ON, Canada}
\email{michelle.delcourt@uwaterloo.ca}

\author{Anush Tserunyan}
\address[Anush Tserunyan]{\normalfont{}Department of Mathematics, University of Illinois at Urbana-Champaign, Urbana, IL, 
	USA}
\email{anush@illinois.edu}

\thanks{Research of the second author is supported in part by EPSRC grant EP/P009913/1 and NSF Graduate Research Fellowship DGE 1144245. Research of the third author is supported in part by NSF Grant DMS-1501036.}

\numberwithin{equation}{section}

\newtheorem{theo}[equation]{Theorem}
\newtheorem{prop}[equation]{Proposition}
\newtheorem{lemma}[equation]{Lemma}
\newtheorem{corl}[equation]{Corollary}

\newtheorem{obs}[equation]{Observation}
\newtheorem{ques}[equation]{Question}

\newcounter{ForClaims}[section]
\newtheorem{claim}{Claim}[ForClaims]

\newcommand*{\myproofname}{Proof}
\newenvironment{claimproof}[1][\myproofname]{\begin{proof}[#1]}{\end{proof}}

\newenvironment{coolproof}[1][Proof]{\begin{proof}[\scshape\upshape#1]}{\end{proof}}

\theoremstyle{definition}
\newtheorem{defn}[equation]{Definition}

\theoremstyle{remark}
\newtheorem*{ques*}{Question}
\newtheorem*{caut*}{Caution}
\newtheorem*{remk*}{Remark}

\renewcommand{\subsetneq}{\varsubsetneq}
\newcommand{\0}{\varnothing}
\newcommand{\set}[1]{{\{#1\}}}
\newcommand{\proj}{\mathrm{proj}}
\newcommand{\dom}{\mathrm{dom}}
\newcommand{\im}{\mathrm{im}}

\newcommand{\N}{\mathbb{N}}
\newcommand{\Z}{\mathbb{Z}}
\newcommand{\F}{\mathbb{F}}
\newcommand{\R}{\mathbb{R}}
\newcommand{\Q}{\mathbb{Q}}
\renewcommand{\c}{^{\mathsf{c}}}
\renewcommand{\U}{\mathscr{U}}

\renewcommand{\epsilon}{\varepsilon}
\renewcommand{\phi}{\varphi}
\renewcommand{\theta}{\vartheta}

\renewcommand{\leq}{\leqslant}
\renewcommand{\geq}{\geqslant}

\renewcommand{\G}{\Gamma}

\newcommand{\B}{\mathscr{B}}

\newcommand{\defeq}{\coloneqq}

\newcommand{\emphd}[1]{\textbf{#1}}

\renewcommand{\d}{\mathfrak{d}}
\renewcommand{\b}{\mathfrak{b}}
\newcommand{\w}{\mathfrak{w}}

\newcommand{\Rat}{\mathcal{R}}
\renewcommand{\L}{\mathcal{L}}

\newcommand{\Comp}{\mathfrak{I}}
\newcommand{\calF}{\mathcal{F}}
\newcommand{\calG}{\mathcal{G}}
\newcommand{\eval}{\mathfrak{e}}
\newcommand{\cl}[1]{\overline{#1}}
\newcommand{\empht}[1]{\textbf{\emph{#1}}}

\newcommand{\bemph}[1]{{\normalfont#1}} 
\newcommand{\ep}[1]{\bemph{(}#1\bemph{)}} 

\renewcommand{\thesubsection}{\arabic{section}.\Alph{subsection}}

\usepackage{etoolbox}
\usepackage{titlesec}

\titleformat{\section}[block]{\scshape\filcenter}{\thesection.}{1ex}{}
\titleformat{\subsection}[block]{\itshape\bfseries\filcenter}{\thesubsection.}{1ex}{}
\titleformat{\subsubsection}[runin]{\scshape}{\thesubsubsection.}{1ex}{}[.]

\titlespacing*{\section}{0pt}{*3}{*1}
\titlespacing*{\subsection}{0pt}{*3}{*1}

\makeatletter
\newcommand{\neutralize}[1]{\expandafter\let\csname c@#1\endcsname\count@}
\makeatother

\newenvironment{theocopy}[1]
{\renewcommand{\thetheo}{\ref{#1}}
	\neutralize{theo}\phantomsection
	\begin{theo}}
	{\end{theo}}

\newenvironment{propcopy}[1]
{\renewcommand{\theprop}{\ref{#1}}
	\neutralize{prop}\phantomsection
	\begin{prop}}
	{\end{prop}}

\renewcommand{\mkbibnamefirst}[1]{\textsc{#1}}
\renewcommand{\mkbibnamelast}[1]{\textsc{#1}}
\renewcommand{\mkbibnameprefix}[1]{\textsc{#1}}
\renewcommand{\mkbibnameaffix}[1]{\textsc{#1}}
\renewcommand{\finentrypunct}{}

\renewbibmacro{in:}{}

\renewbibmacro*{volume+number+eid}{%
	\printfield{volume}%
	\setunit*{\addnbspace}
	\printfield{number}%
	\setunit{\addcomma\space}%
	\printfield{eid}}

\DeclareFieldFormat[article]{volume}{\textbf{#1}\space}
\DeclareFieldFormat[article]{number}{\mkbibparens{#1}}

\DeclareFieldFormat{journaltitle}{#1,}
\DeclareFieldFormat[thesis]{title}{\mkbibemph{#1}\addperiod}
\DeclareFieldFormat[article, unpublished, thesis]{title}{\mkbibemph{#1},}
\DeclareFieldFormat[book]{title}{\mkbibemph{#1}\addperiod}
\DeclareFieldFormat[unpublished]{howpublished}{#1, }

\DeclareFieldFormat{pages}{#1}

\DeclareFieldFormat[article]{series}{Ser.~#1\addcomma}

\renewcommand*{\bibfont}{\small}

\begin{document}
	
	\maketitle
	
	\begin{abstract}
		In this paper we study hypergraphs definable in an algebraically closed field. Our goal is to show, in the spirit of the so-called transference principles in extremal combinatorics, that if a given algebraic hypergraph is ``dense'' in a certain sense, then a generic low-dimensional subset of its vertices induces a subhypergraph that is also ``dense.'' \ep{For technical reasons, we only consider low-dimensional subsets that are parameterized by rational functions.} Our proof approach is inspired by the hypergraph containers method, developed by Balogh, Morris, and Samotij and independently by Saxton and Thomason \ep{although adapting this method to the algebraic setting presents some unique challenges that do not occur when working with finite hypergraphs}. Along the way, we establish a natural generalization of the classical dimension of fibers theorem in algebraic geometry, which is interesting in its own right.
	\end{abstract}
	
	\section{Introduction}
	
	
	An active line of inquiry in combinatorics in recent years has been extending classical results to the so-called \emph{sparse random setting}, where the goal is to show that certain known properties of ``dense'' combinatorial structures are inherited by their randomly chosen ``sparse'' substructures. A typical example of this is the \emph{Sparse Szemer\'edi Theorem}, proved by {Conlon} and {Gowers} \cite{ConlonGowers} and independently by {Schacht} \cite{Schacht}. Say that a finite set $S \subset \N$ is \emphd{$(\epsilon, t)$-Szemer\'edi}, where $\epsilon > 0$ and $t \in \N$, if every subset $A \subseteq S$ with $|A| \geq \epsilon |S|$ contains a nonconstant arithmetic progression of length $t$. The following is a classical theorem of {Szemer\'edi}:
	
	\begin{theo}[{Szemer\'edi \cite{Sz}}]\label{theo:Sz}
		For all $\epsilon > 0$ and $t \in \N$, there exists $n_0 \in \N$ such that for each $n \geq n_0$, the set $\set{1, \ldots, n}$ is $(\epsilon, t)$-Szemer\'edi.
	\end{theo}
	
	\noindent Conlon--Gowers and Schacht proved that the property of being $(\epsilon, t)$-Szemer\'edi is inherited, with high probability, by a fairly sparse random subset of $\set{1, \ldots, n}$:
	
	\begin{theo}[{\empht{Sparse Szemer\'edi}; Conlon--Gowers \cite{ConlonGowers}, Schacht \cite{Schacht}}]\label{theo:sparse_Sz}
		For all $\epsilon > 0$ and $t \in \N$, there exist $n_0 \in \N$ and $C > 0$ such that the following holds:
		
		\smallskip
		
		For each $n \geq n_0$, fix some $p_n \in [Cn^{-1/(t-1)}, 1]$. Let $S_n$ be the random subset of $\set{1, \ldots, n}$ obtained by picking each element $i \in \set{1, \ldots, n}$ independently with probability $p_n$. Then
		\[
			\lim_{n \to \infty} \mathbb{P} \left[\text{$S_n$ is $(\epsilon, t)$-Szemer\'edi}\,\right] \,=\, 1.
		\]
	\end{theo}
	\begin{remk*}
		The lower bound on $p_n$ in Theorem~\ref{theo:sparse_Sz} in sharp, up to the value of $C$.	
	\end{remk*}
	
	\noindent For further examples see, e.g., the survey \cite{ConlonSurvey}.
	
	The remarkable success of this research program is largely due to the development of powerful general techniques for proving random analogs of combinatorial theorems. One of them is the so-called \emph{\ep{hypergraph} containers method}, introduced independently by {Balogh}, {Morris}, and {Samotij} \cite{BMS} and {Saxton} and {Thomason} \cite{ST} and based on the previous work of {Kleitman} and {Winston} \cite{KW} and Sapozhenko \cite{Sapozhenko}.
	A \emphd{\ep{$t$-uniform} hypergraph} on a set $V$ of \emphd{vertices} is a family $E$ of $t$\-/element subsets of $V$, called the \emphd{edges} of $E$. A set $I \subseteq V$ of vertices is \emphd{$E$-independent} if it does not include any edge of $E$ as a subset. The starting point of the containers method is the observation that problems in extremal combinatorics often involve independent sets in specific hypergraphs \ep{for instance, \hyperref[theo:Sz]{Szemer\'edi's theorem} is a statement about independent sets in the $t$-uniform hypergraph on $\set{1, \ldots, n}$ whose edges are the $t$-term arithmetic progressions}. Assuming that the edges of a given hypergraph $E$ are ``well-distributed'' in a certain technical sense, the containers method puts an upper bound on the number of $E$-independent sets, and, furthermore, it does so in a very ``explicit'' manner. Namely, each independent set $I$ gets assigned a \emph{fingerprint} $F$ that, in turn, encodes a \emph{container} $C$, with the property that $F \subseteq I \subseteq C$, meaning that the total number of independent sets is at most the number of distinct fingerprints times the maximum number of independent sets in an individual container.  

	
	In this paper we combine the containers method with another recent trend in combinatorics: establishing versions of extremal results for \ep{hyper}graphs that are definable in model\-/theoretically tame structures; see, e.g., \cite{AlonPachEtc, FoxGromovEtc, Suk, Tao, CS} and the references therein for a sample of related work. We shall specifically focus on the case when the ambient structure is an algebraically closed field \ep{although it would be interesting to know if our results could be interpreted and proved in some more general context}. 
	
	Let $\F$ be an algebraically closed field and let $E$ be a $t$-uniform hypergraph with vertex set $\F^n$. For convenience, we shall view the edges of $E$ as \emph{ordered} tuples of length $t$ rather than simply $t$-element sets, so $E$ is a subset of $(\F^n)^t \cong \F^{tn}$. Thus, it makes sense to ask whether $E$ is a \emph{definable} set \ep{in the sense of first-order logic} in the field structure of $\F$. For instance, if $E$ is the set of all $t$-term arithmetic progressions in $\F^n$, then $E$ {is} definable, since
	\[
		(x_1, \ldots, x_t) \in E \quad \Longleftrightarrow \quad (x_3-x_2 = x_2 - x_1) \wedge \ldots \wedge (x_t - x_{t-1} = x_{t-1} - x_{t-2}) \wedge (x_1 \neq x_2). 
	\]
	
	Given a definable hypergraph $E$ on $\F^n$, we wish to study the properties of the subhypergraph of $E$ induced by a ``sparse random'' \emph{definable} set $X \subseteq \F^n$. It is fairly clear what the word ``sparse'' should mean in this context: Since we are working in an algebraically closed field, there is a notion of \emph{dimension} of definable sets, so ``sparse'' stands for ``low-dimensional.'' It is somewhat less obvious how to interpret the word ``random'' correctly. The approach we take in this paper is to consider only those subsets of $\F^n$ that are \emph{parameterized} by rational functions. Each such subset is described via a finite tuple of elements of $\F$---namely the coefficients of the parameterizing polynomials---and hence can be encoded as a vector in the affine space $\F^N$ for some $N \in \N$. This enables us to talk about a subset $X \subseteq \F^n$ whose parameterization is given by \emph{generic} $x \in \F^N$, and that will be our notion of ``randomness.'' \ep{For details, see \S\ref{subsec:rat_fn}.}
	
	The motivating question can now be stated in reasonably precise terms:
	
	\begin{ques}\label{ques:ques}
		Suppose $E \subseteq (\F^n)^t$ is a definable $t$-uniform hypergraph. What can be said about the subhypergraph $E[X]$ induced by a generic parameterized subvariety $X \subseteq \F^n$ of a given dimension $k$? In particular, what properties of $E$ would guarantee that every definable $E[X]$-independent set $I \subseteq X$ has dimension less than $k$?
	\end{ques}
	
	
	Here and in what follows we use the words ``less'' and ``greater'' to mean ``strictly less'' and ``strictly greater,'' respectively. Our main result is Theorem~\ref{theo:main}, which 
	gives an answer to Question~\ref{ques:ques} that is inspired by the analogous results for finite hypergraphs. Our proof strategy is to adapt the ideas of the containers method and to control the subhypergraph $E[X]$ using a sequence of small ``fingerprints.'' However, while the standard hypergraph containers method involves the so-called \emph{scythe algorithm}---an iterative procedure that runs through the vertices of the hypergraph---we must use a different approach, since our hypergraphs typically have infinitely many vertices. Similar considerations in a different infinitary setting have previously led the current authors together with Henry Towsner \cite{ContFin} to develop a \emph{nonalgorithmic} proof of the containers theorem for finite hypergraphs. Unfortunately, the argument in \cite{ContFin} is still too ``discrete'' and not directly applicable in our current framework. Thankfully, there are other tools available in the algebraic context, most notably irreducibility, that allow us to replace induction over the vertex set with induction on the dimension.
	
	A crucial role in our arguments is played by Theorem~\ref{theo:proj}---a certain generalization of the classical dimension of fibers theorem in algebraic geometry, that is interesting in its own right. We state and prove Theorem~\ref{theo:proj} in Section~\ref{sec:proj}; see \S\ref{subsec:proj_intro} for a motivational discussion. The remainder of the proof of Theorem~\ref{theo:main} proceeds via a sequence of applications of Theorem~\ref{theo:proj} and is presented in Section~\ref{sec:proof}.
	

	\section{Main definitions and results}\label{sec:defns}%
	
	\subsection{Basic notation, terminology, and conventions}\label{subsec:conv}
	
	\subsubsection*{Integers}
	
	We use $\N$ to denote the set of all nonnegative integers. For $n \in \N$, let $[n] \defeq \set{1, \ldots, n}$. 
	By default, the variables $d$, $i$, $j$, $k$, $m$, $n$, $r$, $s$, $t$ range over $\N$.
	
	\subsubsection*{Projections}
	
	For a family of sets $(X_i)_{i \in I}$ and $\0 \neq S \subseteq I$, let
	$
	\proj_S \colon \prod_{i \in I} X_i \to \prod_{i \in S} X_i
	$
	be the projection onto the set $S$ of coordinates. For brevity, given $i \in I$, we write $\proj_i$ instead of $\proj_{\set{i}}$. The Cartesian power $X^t$ of a set $X$ is viewed as a product indexed by $[t]$.
	
	\subsubsection*{Definable sets}
	
	Throughout, we work in a fixed algebraically closed field $\F$. 
	The word ``definable'' always means ``definable in $(\F, +, \cdot)$ with parameters.'' We mostly work in affine spaces $\F^n$, but sometimes we shall also use the projective $n$-space $\mathbb{P}^n$ over $\F$. The only topology we refer to is the Zariski topology \ep{on $\F^n$, $\mathbb{P}^n$, etc.}. The closure of a set $X \subseteq \F^n$ \ep{in the Zariski topology on $\F^n$} is denoted by $\cl{X}$. We say that a definable set $X \subseteq \F^n$ is \emphd{irreducible} if $\cl{X}$ is an irreducible variety. An \emphd{\ep{irreducible} component} of a definable set $X$ is any set of the form $C \cap X$, where $C \subseteq \cl{X}$ is an irreducible component of $\cl{X}$. The \emphd{dimension} $\dim X$ of a definable set $X$ is equal to $\dim \cl{X}$, the Krull dimension of the variety $\cl{X}$ \ep{which coincides with the Morley rank of $X$, see \cite[\S{}6]{Marker}}. 
	By convention, $\dim \0 \defeq -1$. 
	Further algebraic\-/geometric and model-theoretic preliminaries are reviewed in Section~\ref{sec:prelim}.
	
	\subsubsection*{Genericity}
	
	For a definable set $X \subseteq \F^n$ and a property $\mathfrak{P}$ of elements of $\F^n$, we say that a \emphd{generic point $x \in X$ satisfies~$\mathfrak{P}$}, in symbols $\forall^\ast x \in X \, (\mathfrak{P}(x))$, if the set $\mathfrak{P}(X) \defeq \set{x \in X\,:\, \mathfrak{P}(x)}$ is definable and $\cl{\mathfrak{P}(X)} = \cl{X}$. This definition is equivalent to an \ep{apparently stronger} requirement that $\mathfrak{P}(X)$ is a definable set that contains a dense relatively open subset of $\cl{X}$. Another convenient way to phrase this is that generic $x \in X$ satisfies $\mathfrak{P}$ if and only if the set $\mathfrak{P}(X)$ is definable and has codimension $0$ in every irreducible component of $X$. 
	See \S\ref{subsec:quantifiers} for some basic properties of the $\forall^\ast$ quantifier.
		
	\subsection{Rational maps}\label{subsec:rat_fn}
	
	
	By a \emphd{rational map} from $\F^k$ to $\F^n$ we mean an expression of the form
	\begin{equation}\label{eq:rat_fn}
	f(x_1, \ldots, x_k) = \left(\frac{p_1(x_1, \ldots, x_k)}{p_0(x_1, \ldots, x_k)}, \,\ldots,\, \frac{p_n(x_1, \ldots, x_k)}{p_0(x_1, \ldots, x_k)}\right),
	\end{equation}
	where $p_0$, $p_1$, \ldots, $p_n \in \F[x_1, \ldots, x_k]$ and $p_0$ is not identically zero. For convenience, we identify each rational map $f$ as in~\eqref{eq:rat_fn} with the tuple of polynomials $(p_0, p_1, \ldots, p_n)$ rather than with the corresponding partial function $\F^k \rightharpoonup \F^n$ (for instance, multiplying every polynomial in \eqref{eq:rat_fn} by the same element of $\F \setminus \set{0,1}$ produces a different rational map). The set of all rational maps from $\F^k$ to $\F^n$ is denoted by $\Rat(k,n)$.
	
	Given a nonzero polynomial $q \in \F[x_1, \ldots, x_k]$, we write $\Rat_d(k, n; q)$ to indicate for the set of all rational maps $f$ of the form~\eqref{eq:rat_fn} with
	\[p_0 = q \qquad\text{and}\qquad \deg(p_i) \leq d \text{ for all } 1 \leq i \leq n.\]
	In particular, $\Rat_d(k,n;1)$ is the set of all \emphd{polynomial maps} from $\F^k$ to $\F^n$ of degree at most $d$. For brevity, let $\L(k, n) \defeq \Rat_1(k,n;1)$ denote the set of all \emphd{affine maps} from $\F^k$ to $\F^n$. We will also use the shortcut $\L_{<k}(n) \defeq \bigcup_{i < k} \L(i, n)$.
	
	A rational map in $\Rat_d(k,n;q)$ can be described by a tuple of
	$
		 N \defeq n {k+d \choose d}
	$
	elements of $\F$, namely the coefficients of the corresponding polynomials $p_1$, \ldots, $p_n$ in \eqref{eq:rat_fn}.
	As mentioned in the introduction, we use this observation to identify $\Rat_d(k,n;q)$ with the space $\F^{N}$, so it makes sense to talk about definable subsets of $\Rat_d(k,n;q)$ as well as the properties of a {generic} element of $\Rat_d(k,n;q)$.
	
	\subsection{Definable independent sets}
	
	\noindent Let $t \geq 1$ and $E \subseteq (\F^n)^t$. We say that a set $I \subseteq \F^n$ is \emphd{$E$-independent} if $E \cap I^t = \0$, i.e., there are no $x_1$, \ldots, $x_t \in I$ with $(x_1, \ldots, x_t) \in E$.
	
	A common feature of the techniques for proving results in the sparse random setting is their reliance on the corresponding theorems in the \emph{dense} case. For instance, all the known proofs of Theorem~\ref{theo:sparse_Sz} treat \hyperref[theo:Sz]{Szemer\'edi's theorem}\footnote{Or, more precisely, its robust version, originally due to Varnavides \cite{Varn} in the case of $3$-term progressions.} as a black box. Proving a result in the dense case can often be a challenging task in its own right; and indeed, the principal obstacles to resolving several open problems in the theory of sparse random structures lie in obtaining sufficiently strong ``supersaturation'' bounds in the dense setting \ep{see, e.g., \cite{BMS_survey, supersat}}. However, the situation simplifies dramatically in the algebraic context, thanks to the following proposition, which provides a convenient criterion for when a definable hypergraph $E$ admits a ``large'' definable independent set:
	
	\begin{prop}\label{prop:dense}
		Let $t \geq 1$ and let $E \subseteq (\F^n)^t$ be a definable set. The following statements are equivalent:
		\begin{enumerate}[label=\ep{\normalfont\roman*}]
			\item\label{item:ind} the dimension of every definable $E$-independent set $I \subseteq \F^n$ is less than $n$;
			
			\item\label{item:proj} $E$ has an irreducible component $H \subseteq E$ such that for all $i \in [t]$, we have $\dim \proj_i H = n$.
		\end{enumerate}
	\end{prop}
	\begin{coolproof}
		First, consider the case when $E$ is irreducible. To prove \ref{item:ind} $\Longrightarrow$ \ref{item:proj}, note that if $\dim \proj_i E < n$ for some $i \in [t]$, then $\F^n \setminus \proj_i E$ is a definable $E$-independent set of dimension $n$. Towards the other implication, assume \ref{item:proj} and let $I \subseteq \F^n$ be a definable set of dimension $n$. For each $i \in [t]$, let \[E_i \defeq E \cap \proj_i^{-1}(I).\] Since, by assumption, $\dim \proj_i E = \dim I = n$ and $E$ is irreducible, the dimension of fibers theorem \ep{specifically, Corollary~\ref{corl:fubini}} yields $\dim E_i = \dim E$. The irreducibility of $E$ then gives \[\dim(E \cap I^t) = \dim (E_1 \cap \ldots \cap E_t) = \dim(E) \geq 0,\]	and hence $I$ is not $E$-independent.
		
		Now let $E$ be arbitrary. Since we have already verified the proposition for irreducible sets, it remains to show that if every definable \emph{$E$-independent} set is of dimension less than $n$, then $E$ has a component $H$ such that every definable \emph{$H$-independent} set is of dimension less than $n$ as well. To that end, suppose that the components of $E$ are $H_1$, \ldots, $H_k$ and, for each $j \in [k]$, let $I_j \subseteq \F^n$ be a definable $H_j$-independent set of dimension $n$. Since $E = H_1 \cup \ldots \cup H_k$, we conclude that $I_1 \cap \ldots \cap I_k$ is a definable $E$-independent set of dimension $n$, which finishes the proof.
	\end{coolproof}
	
	An immediate consequence of Proposition~\ref{prop:dense} is, for example, the fact that for any $t \geq 3$, there is no $n$-dimensional definable set $I \subseteq \F^n$ without a $t$-term arithmetic progression, since the set
	\[
		E \defeq \set{(x_1, \ldots, x_t) \in (\F^n)^t \,:\, x_2 - x_1 = \cdots = x_t - x_{t-1} \text{ and } x_1 \neq x_2}
	\]
	is irreducible and satisfies $\proj_i E = \F^n$ for all $i \in [t]$. It might perhaps seem strange that this algebraic version of \hyperref[theo:Sz]{Szemer\'edi's theorem} is almost trivial, while \hyperref[theo:Sz]{Szemer\'edi's theorem} itself is an extremely deep result. The explanation is simple. Any {definable} set $I \subseteq \F^n$ must fulfill one of the two alternatives: either $\dim I < n$, 
	or else, $\dim (\F^n \setminus I) < n$. This is analogous, in the finite setting, to $I$ either having \emph{density $0$} or \emph{density $1$}. On the other hand, the statement of \hyperref[theo:Sz]{Szemer\'edi's theorem} is only difficult for sets whose density is \emph{small but positive}---and they do not exist in the algebraic setting.
	
	What makes Proposition~\ref{prop:dense} particularly useful is that, while property~\ref{item:ind} is of primary interest to us, it is not apparently first-order, due to the quantification over all definable sets $I \subseteq \F^n$; on the other hand, property~\ref{item:proj} is more ``concrete,'' and indeed, it is definable in families \ep{see Corollary~\ref{corl:large_ind_set_def}}. Notice however that, to be able to verify \ref{item:proj}, we must have good control over the individual irreducible components of $E$, as it is not enough to simply know that $\dim \proj_i E = n$ for all $i \in [t]$. For instance, it may happen that $E$ is expressible as a union $E = E_1 \cup E_2$ of two definable sets $E_1$, $E_2 \subseteq (\F^n)^t$ with $\dim\proj_1 E_1 = \dim \proj_2 E_2 = n$ but $\dim \proj_2 E_1$, $\dim \proj_1 E_2 < n$, in which case property \ref{item:proj}, and hence also \ref{item:ind}, fails. \ep{The simplest example of this situation is the set $E = (\F^n \times \set{0}) \cup (\set{0} \times \F^n) \subseteq (\F^n)^2$.} This issue will be a source of some important technical complications in our arguments.

	\subsection{The main result}
	
	Let $E \subseteq (\F^n)^t$ be a definable set. Given a definable map $f \colon \F^k \rightharpoonup \F^n$, we wish to consider the subhypergraph of $E$ induced by the subset $f(\F^k) \subseteq \F^n$ parameterized by $f$. Since it will be more convenient to work directly in the parameter space $\F^k$, we define the \emphd{subhypergraph of $E$ induced by $f$} to be the set $E[f] \subseteq (\F^k)^t$ given by
	\[
	(y_1, \ldots, y_t) \in E[f] \quad \vcentcolon\Longleftrightarrow \quad (f(y_1), \ldots, f(y_t)) \in E.
	\]
	A peculiar feature of this definition is that it makes sense even when $k > n$; in other words, the dimension of the vertex set of $E[f]$ can exceed that of the vertex set of $E$.
	
	We say that a definable set $E \subseteq (\F^n)^t$ is \emphd{injective} if for all $(x_1, \ldots, x_t) \in E$, the elements $x_1$, \ldots, $x_t$ are pairwise distinct. When $E$ is thought of as a hypergraph on $\F^n$, the injectivity of $E$ means that it is ``truly'' $t$-uniform, i.e.,  every edge of $E$ contains precisely $t$ distinct vertices.
	
	The next definition is particularly important. We say that a definable set $E \subseteq (\F^n)^t$ is \emphd{$r$-almost dense} if for all subsets $\0 \neq S \subseteq [t]$, we have
	\begin{equation}\label{eq:r-dense}
			\dim \proj_S E \,\geq\, |S|n - (|S|-1)r.
	\end{equation}	
	Observe the following chain of implications:
	\[
		\cl{E} = (\F^n)^t \ \Longleftrightarrow\ \text{$E$ is $0$-almost dense} \ \Longrightarrow\ \text{$E$ is $1$-almost dense} \ \Longrightarrow\ \text{$E$ is $2$-almost dense} \ \Longrightarrow\ \ldots
	\]
	The notion of $r$-almost density for subsets $E \subseteq (\F^n)^t$ is mostly interesting for $r \leq n$, as if $r \geq n$, then \begin{equation}\label{eq:dp}
	E \subseteq (\F^n)^t \text{ is $r$-almost dense} \quad \Longleftrightarrow \quad \dim \proj_i E = n \text{ for all } i \in [t].
	\end{equation}
	We emphasize that in \eqref{eq:r-dense}, $n$ is the dimension of the vertex set of $E$; for instance, when we say that a definable set $E \subseteq (\F^k)^t$ is $r$-almost dense, then $n$ must be replaced by $k$.
	
	\begin{obs}\label{obs:ind}
		Let $t \geq 1$ and let $E \subseteq (\F^n)^t$ be an irreducible definable set. If $E$ is $n$-almost dense, then the dimension of every definable $E$-independent set $I \subseteq \F^n$ is less than $n$.
	\end{obs}
	\begin{coolproof}
		Immediate from \eqref{eq:dp} and Proposition~\ref{prop:dense}.
	\end{coolproof}
	
	The significance of the notion of almost density is demonstrated by the following proposition, whose proof is deferred to \S\ref{subsec:lower_bound}:
	
	\begin{prop}\label{prop:lower_bound}
		Let $t \geq 1$ and let $E \subseteq (\F^n)^t$ be an injective definable set. If there exist $d \geq t-1$ and a nonzero polynomial $q \in \F[x_1, \ldots, x_k]$ such that, for generic $f \in \Rat_d(k,n;q)$, every definable $E[f]$-independent set $I \subseteq \F^k$ has dimension less than $k$, then $E$ has a $k$-almost dense irreducible component.
	\end{prop}

	Now we are ready to state the main result of this paper, answering Question~\ref{ques:ques}:
	
	\begin{theo}\label{theo:main}
		Let $t \geq 1$ and let $E \subseteq (\F^n)^t$ be an $r$-almost dense irreducible definable set. Fix $d \geq t-1$ and a nonzero polynomial $q \in \F[x_1, \ldots, x_k]$. If $k \geq r + 1$, then, for generic $f \in \Rat_d(k,n;q)$, the following holds:
		
		\smallskip
		
		Every definable $E[f]$-independent set $I \subseteq \F^k$ has dimension less than $k$. Furthermore, if $E$ is injective, then every irreducible component of $E[f]$ is $r$-almost dense.
	\end{theo}
	
	We finish this section with a few remarks about the statement of Theorem~\ref{theo:main}.
	
	
	
	\subsubsection*{Lower bound on $d$}
	
	Theorem~\ref{theo:main} requires $d$ to grow with $t$, the uniformity of the hypergraph $E$. Informally, one could say that the maps given by polynomials of low degree are not ``random enough'' for the conclusion of Theorem~\ref{theo:main} to hold. In fact, the lower bound $d \geq t-1$ is best possible, as the following construction shows. Take any $t \geq 2$, $d \leq t-2$, and $n > t$. Define
	\[
		H \defeq \set{(x_1, \ldots, x_t) \in (\F^n)^t \,:\, \exists f \in \Rat_d(1, n;1) \, (x_1, \ldots, x_t \in f(\F))}.
	\]
	Set $E \defeq (\F^n)^t \setminus H$. Using that $d \leq t-2$ and $n > t$, we obtain
	\[
	\dim H \,\leq\, \dim (\Rat_d(1,n;1) \times \F^t) \,=\, n(d+1) + t \,\leq\, n(t-1) + t\,<\, nt \,=\, \dim (\F^n)^t,
	\]
	and hence $\cl{E} = (\F^n)^t$. But, on the other hand, $E[f] = \0$ {for all} $f \in \Rat_d(1, n; 1)$ by definition.
	
	\subsubsection*{Lower bound on $k$}
	
	If, in the setting of Theorem~\ref{theo:main}, we let $r$ be the smallest integer such that $E$ is $r$-almost dense, then it follows from Proposition~\ref{prop:lower_bound} that $k$ must be at least $r$. Actually, the conclusion of Theorem~\ref{theo:main} can fail even when $k = r$ \ep{and hence the lower bound $k \geq r + 1$ is best possible}. To see this, take any $k < n$ and define a set $E \subseteq (\F^n)^2$ as follows. Pick an arbitrary $(n-k)$-dimensional linear subspace $V \subseteq \F^n$ and put
	\[
		(x, y) \in E \quad \vcentcolon\Longleftrightarrow \quad x \neq y \text{ and } x-y \in V.
	\]
	Since $n-k > 0$, we have $\proj_1 E = \proj_2 E = \F^n$ and $\dim E =  2n - k$, so $E$ is $k$-almost dense. However, for a generic affine map $\ell \in \L(k, n)$, the set $\ell(\F^k)$ is a $k$-dimensional affine subspace of $\F^n$ that intersects every translate of $V$ in precisely one point---and hence $E[\ell]=\0$.
	
	\subsubsection*{The codegree conditions}	
	
	We now briefly comment on the relationship between the statement of Theorem~\ref{theo:main} and the sparse random results in the finite setting.
	
	In order to establish a transference principle for a {finite} $t$-uniform hypergraph $E$, it is typical to assume that $E$ has ``many'' edges, and that the edges of $E$ are somehow ``well-distributed'' over the vertex set. For example, in the containers method, these assumptions take the form of the \emph{codegree conditions}: The degree of each vertex of $E$ is required to be close to the average, and similar restrictions are put on the codegrees of all sets of fewer than $t$ vertices.
	
	In Theorem~\ref{theo:main}, the part that forces the edges of $E$ to be ``well-distributed'' is the assumption of \emph{irreducibility}. For instance, it follows from the dimension of fibers theorem \ep{see \S\ref{subsec:fib_dim}} that if a definable set $E \subseteq (\F^n)^t$ is irreducible and satisfies $\dim \proj_i E = n$ for all $i \in [t]$, then the ``degree'' of a generic vertex $x \in \F^n$---i.e., the dimension of the set of all tuples $(x_1, \ldots, x_t) \in E$ such that $x \in \set{x_1, \ldots, x_t}$---is equal to the ``average'' value $\dim E - n$. Since irreducibility takes care of the ``uniform distribution'' of the edges, there is no need to explicitly bound the codegrees in the statement of Theorem~\ref{theo:main}, and the only numerical assumption left is that $E$ has ``many'' edges---specifically, it must be $r$-almost dense.
	
	Nevertheless, there are still some close parallels between the statement of Theorem~\ref{theo:main} and, say, that of \cite[Proposition~3.1]{BMS}. Indeed, let $E$ be a finite $t$-uniform hypergraph on a set $X$ and imagine that we wish to apply \cite[Proposition~3.1]{BMS} to $E$; in particular, let $p$ be a value between $0$ and $1$. The conclusion of \cite[Proposition~3.1]{BMS} is only interesting for independent sets of size at least $\Omega(p |X|)$. Set $n \defeq \log |X|$ and $r \defeq \log (p|X|)$. For each $s \in [t]$, let $m_s$ be the logarithm of the number of $s$-element subsets of $X$ that are contained in at least one edge of $E$. The codegree conditions of \cite[Proposition~3.1]{BMS} then imply
	\[
		m_s \,\geq\,sn -  (s - 1) r + O(1),
	\]
	which should be compared to \eqref{eq:r-dense}.
	
	\section{ Preliminaries}\label{sec:prelim}
	
	\subsection{Definability in algebraically closed fields}\label{subsec:model_theory}
	
	The following fundamental fact and its immediate consequences will be used without mention:
	
	\begin{theo}[{\empht{Quantifier elimination}~\cite[Theorem~3.2.2]{Marker}}]\label{theo:QE}
		Every definable set $X \subseteq \F^n$ is constructible, i.e., it is a finite Boolean combination of closed sets.
	\end{theo}
	
	For $A \subseteq X \times Y$ and $x \in X$, we let $A_x$ denote the \emphd{fiber of $A$ over $x$}, i.e., the set
	\[
		A_x \defeq \set{y \in Y \,:\, (x, y) \in A}.
	\]
	The next theorem follows from \cite[\S{}I.8, Corollary~3]{Mum} and is a special case of the fact that Morley rank is definable in strongly minimal theories \cite[Lemma 6.2.20]{Marker}:
	
	\begin{lemma}[{\empht{Definability of dimension} \cite[\S{}I.8, Corollary~3]{Mum}}]\label{lemma:dim}
		Let $A \subseteq \F^n \times \F^m$ be a definable set and let $k \in \N$. Then the set $\set{x \in \F^n \,:\, \dim A_x = k}$ is definable.
	\end{lemma}
	
	We shall also require certain more advanced facts concerning definability in algebraically closed fields, all of which are well-known and can be easily derived from general results obtained by van den Dries and Schmidt in their foundational paper \cite{SLou}. For additional information and alternative presentation, we also refer the reader to Chapter IV of van den Dries's thesis \cite{Lou} and to Chapter 10 of Johnson's thesis \cite{Will}.
	
	\begin{lemma}[{\empht{Fiberwise closure} \cite[Theorem~10.2.1(1)]{Will}}]\label{lemma:closure}
		Let $A \subseteq \F^n \times \F^m$ be a definable set. Then the set $\set{(x,y) \in \F^n \times \F^m \,:\, y \in \cl{A_x}}$ is definable.
	\end{lemma}
	
	\begin{lemma}[{\empht{Irreducible fibers} \cite[Theorem~10.2.1(2)]{Will}}]\label{lemma:irred_fibers}
		Let $A \subseteq \F^n \times \F^m$ be a definable set. Then the set $\set{x \in \F^n \,:\, \text{$A_x$ is irreducible}}$ is definable.
	\end{lemma}
	
	
	The next definability result will play a particularly important role in the sequel. For a definable set $X \subseteq \F^n$, let $\Comp(X) \subseteq X \times X$ denote the relation given by
	\[
	(x_1, x_2) \in \Comp(X) \,\vcentcolon\Longleftrightarrow\, \text{there is an irreducible component $I$ of $X$ such that $x_1$, $x_2 \in I$}.
	\]
	Note that $\Comp(X)$ is ``almost'' an equivalence relation: it is reflexive and symmetric, but not necessarily transitive, as distinct irreducible components of $X$ need not be disjoint.
	
	\begin{lemma}[\empht{Definability of components}]\label{lemma:comps_are_definable}
		Let $A \subseteq \F^n \times \F^m$ be a definable set. Then the set
		\[
			\mathfrak{X} \defeq \set{(x, y_1, y_2) \in \F^n \times \F^m \times \F^m \,:\, (y_1, y_2) \in \Comp(A_x)}
		\]
		is definable as well.
	\end{lemma}
	\begin{coolproof}
		Lemma~\ref{lemma:comps_are_definable} is a straightforward consequence of \cite[Theorem~2.10]{SLou} and other results of that paper. For completeness, we include a direct derivation from Lemmas~\ref{lemma:closure} and \ref{lemma:irred_fibers} here. This is the only place in the paper where we invoke nontrivial model-theoretic machinery. Replacing $\F$ by an elementary extension if necessary, we may arrange that $\F$ is $\aleph_1$-saturated.\footnote{Actually, \emph{every} {uncountable} algebraically closed field is saturated, since the theory of algebraically closed fields of any fixed characteristic is uncountably categorical, hence $\omega$-stable; see \cite[Theorem 6.5.4]{Marker}.}\footnote{Here we use the fact that if $\mathbb{K} \supseteq \F$ is an algebraically closed field extending $\F$ and $X \subseteq \mathbb{K}^n$ is a set definable in $\mathbb{K}$, then $X \cap \F^n$ is definable in $\F$. This can be verified using quantifier elimination \ep{i.e., Theorem~\ref{theo:QE}} and the fact that if $X \subseteq \mathbb{K}^n$ is closed in $\mathbb{K}^n$, then $X \cap \F^n$ is closed in $\F^n$. An analogous result holds more generally in arbitrary stable theories and is a consequence of definability of types; see \cite[Exercise 1.29]{Pillay}.} 
		For a definable set $Y \subseteq \F^m$, let $\Comp_{k, d}(Y)$ denote the set of all pairs $(y_1, y_2) \in Y \times Y$ such that $y_1$ and $y_2$ belong to an irreducible component of $\cl{Y}$ cut out from $\F^m$ by at most $k$ polynomials of degree at most $d$. 
		Clearly, $\Comp(Y) = \bigcup_{k, d} \Comp_{k, d}(Y)$, and, since a closed set $Z \subseteq \F^m$ is an irreducible component of $\cl{Y}$ if and only if $Z$ is irreducible, $Z \subseteq \cl{Y}$, and $Z \not \subseteq \cl{Y \setminus Z}$, by Lemmas~\ref{lemma:closure} and \ref{lemma:irred_fibers}, the set
		\[
			\mathfrak{X}_{k,d} \defeq \set{(x, y_1, y_2) \in \F^n \times \F^m \times \F^m \,:\, (y_1, y_2) \in \Comp_{k,d}(A_x)}
		\]
		is definable. Hence, $\mathfrak{X} = \bigcup_{k,d} \mathfrak{X}_{k,d}$ is a countable \emph{union} of definable sets.
		On the other hand, let $\Comp^{k,d}(Y)$ be the set of all pairs $(y_1, y_2) \in Y \times Y$ such that every closed subset $Z \subseteq \F^m$ cut out by at most $k$ polynomials of degree at most $d$ satisfies 
		$
		(y_1 \not\in \cl{Y \setminus Z} \, \Longrightarrow \, y_2 \in Z)	
		$.
		Then $\Comp(Y) = \bigcap_{k,d} \Comp^{k,d}(Y)$, and Lemma~\ref{lemma:closure} yields that the set
		\[
			\mathfrak{X}^{k,d} \defeq \set{(x, y_1, y_2) \in \F^n \times \F^m \times \F^m \,:\, (y_1, y_2) \in \Comp^{k,d}(A_x)}
		\]
		is definable, so $\mathfrak{X} = \bigcap_{k,d} \mathfrak{X}^{k,d}$ is a countable \emph{intersection} of definable sets. Since $\F$ is $\aleph_1$-saturated, a set that is both a countable union and a countable intersection of definable sets must itself be definable.
	\end{coolproof}
	 
	Lemma~\ref{lemma:comps_are_definable} allows using quantification over irreducible components when forming a definable set. Below we present, as an illustration, a typical example of how Lemma~\ref{lemma:comps_are_definable} can be applied: 
	
	\begin{corl}\label{corl:large_ind_set_def}
		Let $E \subseteq \F^n \times (\F^m)^t$ be a definable set. Then the following set is definable:
		\[
			\set{x \in \F^n \,:\, \text{the dimension of every definable $E_x$-independent set $I \subseteq \F^m$ is less than $m$}}.
		\]
	\end{corl}
	\begin{coolproof}
		By Proposition~\ref{prop:dense}, the set in question coincides with the set of all $x \in \F^n$ such that $E_x$ has an irreducible component $H$ with $\dim \proj_i H = m$ for all $i \in [t]$. The existential quantification over the components of $E_x$ can be turned into a quantification over the points of $E_x$ as follows: There is an element $e \in E_x$ such that:
		\begin{enumerate}[label=\ep{\emph{\alph*}}]
			\item\label{item:a} for all $h_1$, $h_2 \in E_x$, if $(e, h_1) \in \Comp(E_x)$ and $(e,h_2) \in \Comp(E_x)$, then $(h_1, h_2) \in \Comp(E_x)$;
			
			\item\label{item:b} for all $i \in [t]$, $\dim \proj_i \set{h \in E_x \,:\, (e, h) \in \Comp(E_x)} = m$.
		\end{enumerate}
		It remains to apply Lemmas~\ref{lemma:comps_are_definable} and \ref{lemma:dim}.
	\end{coolproof}
	
	Throughout the rest of this paper, arguments that are similar to and just as straightforward as the proof of Corollary~\ref{corl:large_ind_set_def} are omitted.
	
	\subsection{The dimension of fibers theorem and its consequences}\label{subsec:fib_dim}
	
	A central role in the sequel is played by the so-called \emph{dimension of fibers theorem}. It is a fundamental result, with many versions and generalizations that fall far beyond the scope of this article. We only give here the statements that will be explicitly used later on; for a more thorough discussion, see, e.g., \cite[\S{}I.8]{Mum} and \cite[\S{}11.4]{Vakil}. 
	

	\begin{theo}[{\empht{Dimension of fibers}, ess.\ \cite[\S{}I.8, Corollary 1]{Mum}}]\label{theo:fib_dim}
		Let $A \subseteq \F^n \times \F^m$ be an irreducible definable set. For generic $x \in \proj_1 A$, the dimension of every component of the set $A_x$ is $\dim A - \dim \proj_1 A$.
	\end{theo}
	
	Since the above theorem is usually stated for varieties rather than definable sets, we include here a derivation of the general case of Theorem~\ref{theo:fib_dim} from the case when $A$ is a closed set: 
	
	\begin{coolproof}[Derivation of Theorem~\ref{theo:fib_dim} from the closed set case]
		Let $B \defeq \cl{A} \setminus A$ and note that $\dim B < \dim A = \dim \cl{A}$.
		Since we assume that Theorem~\ref{theo:fib_dim} holds for $\cl{A}$, it suffices to argue that
		\begin{equation}\label{eq:closure}
			\forall^\ast x \in \proj_1 \cl{A} \qquad \cl{A_x} = (\cl{A})_x.
		\end{equation}
		Suppose that \eqref{eq:closure} fails. Since the set $\proj_1 \cl{A}$ is irreducible, we then have
		\[
			\forall^\ast x \in \proj_1 \cl{A} \qquad \cl{A_x} \subsetneq (\cl{A})_x.
		\]
		Hence, for generic $x \in \proj_1 \cl{A}$, at least one component of the fiber $(\cl{A})_x$ is entirely contained in $\cl{B_x} \subseteq (\cl{B})_x$. In particular, $\dim \proj_1 \cl{A} = \dim \proj_1 \cl{B} \eqqcolon k$. However, Theorem~\ref{theo:fib_dim}, applied to the closed sets $\cl{A}$ and $\cl{B}$, yields that, for generic $x \in \proj_1 \cl{A}$, the dimension of every component of $(\cl{A})_x$ is $\dim A - k$, while $\dim (\cl{B})_x = \dim B - k < \dim A - k$.
		This contradiction completes the proof of \eqref{eq:closure}.
	\end{coolproof}
	
	Claim \eqref{eq:closure} is useful enough to be stated as a separate corollary:
	
	\begin{corl}\label{corl:fibers_of_the_closure}
		Let $A \subseteq \F^n \times \F^m$ be a definable set. Then, for generic $x \in \proj_1 \cl{A}$, we have $\cl{A_x} = (\cl{A})_x$.
	\end{corl}
	
	The next result is a well-known and easy consequence of Theorem~\ref{theo:fib_dim}: 
	
	\begin{corl}[\empht{Fubini for dimension}]\label{corl:fubini}
		Let $A \subseteq \F^n \times \F^m$ be a nonempty definable set. Then
		\[
			\dim A \,=\, \max_{k \in \N} \,\left(k \,+\, \dim \set{x \in \proj_1 A \,:\, \dim A_x = k}\right).
		\] 
	\end{corl}
	
	The following immediate consequence of Corollary~\ref{corl:fibers_of_the_closure} will be used repeatedly: 
	
	\begin{corl}\label{corl:proj_cl}
		Let $A \subseteq \F^n \times \F^m$ be a definable set. Then $\cl{\proj_1 A} = \cl{\proj_1 \cl{A}}$.
	\end{corl}
	
	The following statement will play a crucial role in the later stages of the proof of Theorem~\ref{theo:main}:
	
	\begin{corl}[\empht{Generic indecomposability of fibers}]\label{corl:split}
		Let $A \subseteq \F^n \times \F^m$ be a nonempty irreducible definable set and let $B \subseteq A$ be a definable subset such that for generic $x \in \proj_1 A$, the set $B_x$ is a union of irreducible components of $A_x$. Then
		\[
			\text{either} \qquad \forall^\ast x \in \proj_1 A \, (B_x = \0)\text{,} \qquad \text{or} \qquad \forall^\ast x \in \proj_1 A \, (B_x = A_x).
		\]
	\end{corl}
	\begin{coolproof}
		Suppose that, for generic $x \in \proj_1 A$, the set $B_x$ is nonempty. Then, by Theorem~\ref{theo:fib_dim}, we must have $\forall^\ast x \in \proj_1 A \, (\dim B_x = \dim A - \dim \proj_1 A)$, and therefore, \[\dim B \,\geq\, \dim \proj_1 A + (\dim A - \dim \proj_1 A) \,=\, \dim A.\] Since $A$ is irreducible, this yields $\cl{B} = \cl{A}$, and hence $\cl{B_x} = \cl{A_x}$ for generic $x \in \proj_1 A$. But if $B_x$ is a union of components of $A_x$, then $\cl{B_x} = \cl{A_x}$ is equivalent to $B_x = A_x$, and we are done.
	\end{coolproof}
	
	
	\subsection{An irreducibility criterion}
	
	To verify that certain sets appearing in the proof of Theorem~\ref{theo:main} are irreducible, we will need the following fact \ep{it is the only statement in this paper that requires leaving the realm of affine spaces}:
	
	\begin{lemma}[{ess.\ \cite[Exercise~11.4.C]{Vakil}}]\label{lemma:same_fibers}
		Let $A \subseteq \F^n \times \F^m$ be a nonempty definable set and let $A^\ast$ be the closure of $A$ in $\F^n \times \mathbb{P}^m$. Suppose that $\proj_1 A$ is irreducible and all the fibers $A^\ast_x \subseteq \mathbb{P}^m$ with $x \in \proj_1 A$ are irreducible and of the same dimension. Then $A$ is irreducible.
	\end{lemma}
	\begin{coolproof}
		Let $X \defeq \proj_1 A$ and let $d$ denote the common dimension of the fibers of $A^\ast$ over the points in $X$. Then we must have $d = \dim A - \dim X$. Let $B$ be any irreducible component of $A^\ast$ such that $\dim \proj_1 B = \dim X$. Since $\mathbb{P}^m$ is a complete variety \cite[\S{}I.9, Theorem~1]{Mum}, $\proj_1 B$ is closed in $\F^n$, and, since $X$ is irreducible, $X \subseteq \proj_1 B$. By \cite[\S{}I.8, Theorem 2]{Mum}, we obtain that for all $x \in X$,
		\[
			\dim B_x \,\geq\, \dim B - \dim \proj_1 B \,=\, \dim A - \dim X \,=\, d,
		\]
		and thus $\dim B_x = \dim A_x^\ast$. Since $A_x^\ast$ is irreducible, this implies $B_x = A_x^\ast$ and hence $B = A^\ast$. 
	\end{coolproof}
	
	
	\begin{corl}\label{corl:linear_irred}
		Suppose that $p_1$, \ldots, $p_k \in \F[x_1, \ldots, x_n, y_1, \ldots, y_m]$ are polynomials that are linear in the set of variables $y_1$, \ldots, $y_m$ and let $Z \subseteq \F^n \times \F^m$ be the set of common zeros of $p_1$, \ldots, $p_k$. If $X \subseteq \F^n$ is an irreducible definable set such that $\dim Z_x$ is the same for all $x \in X$, then $Z \cap (X \times \F^m)$ is irreducible.
	\end{corl}
	\begin{coolproof}
		If $Z_x = \0$ for all $x \in X$, then $Z \cap (X \times \F^m) = \0$ and there is nothing to prove. Otherwise, let $Z^\ast$ be the closure of $Z$ in $\F^n \times \mathbb{P}^m$. For each $x \in \F^n$, the set $Z^\ast_x$ is cut out from $\mathbb{P}^m$ by a system of homogeneous linear equations. Hence, the fiber $Z^\ast_x$ is irreducible and if $Z_x \neq \0$, then $Z_x$ is dense in $Z^\ast_x$. The desired conclusion now follows by Lemma~\ref{lemma:same_fibers}.
	\end{coolproof}
	
	\subsection{Fubini-like properties of the generic quantifier}\label{subsec:quantifiers}
	
	Now we state and prove some results that allow changing the order of multiple $\forall^\ast$ quantifiers, in the spirit of the theorems of Fubini in measure theory and Kuratowski--Ulam in general topology.
	
	\begin{theo}[\empht{Fubini for $\forall^\ast$}]\label{theo:Fubini}
		Let $A \subseteq \F^n \times \F^m$ be an irreducible definable set and let $X \defeq \proj_1 A$. Then, for any definable set $B \subseteq A$, the following statements are equivalent:
		\begin{enumerate}[label=\ep{\normalfont\roman*}]
			\item\label{item:1} $\forall^\ast x \in X \quad \forall^\ast y \in A_x \quad ((x,y) \in B)$;
			
			\item\label{item:2} $\forall^\ast x \in X \quad \dim B_x = \dim A_x$;
			
			\item\label{item:3} $\forall^\ast (x,y) \in A\, ((x,y) \in B)$, i.e., $\dim B = \dim A$.
		\end{enumerate}
	\end{theo}
	\begin{coolproof}
		Implication \ref{item:1} $\Longrightarrow$ \ref{item:2} is clear. For \ref{item:2} $\Longrightarrow$ \ref{item:3}, we use the \hyperref[theo:fib_dim]{dimension of fibers theorem} to observe that $\forall^\ast x \in X \, (\dim A_x = \dim A - \dim X)$, so if \ref{item:2} holds, then $\dim B \geq \dim X + (\dim A - \dim X) = \dim A$, as desired. Finally, to prove \ref{item:3} $\Longrightarrow$ \ref{item:1}, consider the set $C \subseteq A$ such that for each $x \in X$, the fiber $C_x$ is the union of all the irreducible components of $A_x$ in which $B_x$ is dense. It follows from the results of \S\ref{subsec:model_theory} that $C$ is definable, and, by Corollary~\ref{corl:split}, we either have $\forall^\ast x \in X \, (C_x = A_x)$ or $\forall^\ast x \in X \, (C_x = \0)$. In the former case, \ref{item:1} holds \ep{and we are done}, so assume that $\forall^\ast x \in X \, (C_x = \0)$, i.e., $\forall^\ast x \in X \, \forall^\ast y \in A_x \, ((x, y) \not \in B)$. Applying \ref{item:1} $\Longrightarrow$ \ref{item:3} with $A \setminus B$ in place of $B$ yields $\dim (A \setminus B) = \dim A$, and hence $\dim B < \dim A$, as desired.
	\end{coolproof}
	
	\begin{corl}\label{corl:FFubini}
		Let $X \subseteq \F^n$, $Y \subseteq \F^m$, and $A \subseteq X \times Y$ be definable sets. Then the following statements are equivalent:
		\begin{enumerate}[label=\ep{\normalfont\roman*}]
			\item\label{item:X_Y} $\forall^\ast x \in X \quad \forall^\ast y \in Y \quad ((x, y) \in A)$;
			
			\item\label{item:Y_X} $\forall^\ast y \in Y \quad \forall^\ast x \in X \quad ((x, y) \in A)$;
			
			\item\label{item:X_times_Y} $\forall^\ast (x, y) \in X \times Y \quad ((x, y) \in A)$.
		\end{enumerate}
	\end{corl}
	\begin{coolproof}
		Since the components of $X \times Y$ are precisely the products of the components of $X$ and $Y$, we may assume that $X$ and $Y$ are irreducible. An application of Theorem~\ref{theo:Fubini} completes the proof.
	\end{coolproof}
	
	In the next corollary, it is instructive to think of $A \subseteq X \times X$ as being an equivalence relation on $X$ \ep{and this is the context in which this corollary will be used later on}.
	
	\begin{corl}\label{corl:equivalence}
		Let $X \subseteq \F^n$ be an irreducible definable set and let $A \subseteq X \times X$ be an irreducible definable set such that $\proj_1 A = \proj_2 A = X$. Then, for any definable set $Y \subseteq X$, the following statements are equivalent:
		\begin{enumerate}[label=\ep{\normalfont\roman*}]
			\item\label{item:E} $\forall^\ast x \in X \quad \forall^\ast x' \in A_x \quad (x' \in Y)$;
			\item\label{item:X} $\forall^\ast x \in X \, (x \in Y)$, i.e., $\dim Y = \dim X$.
		\end{enumerate}
	\end{corl}
	\begin{coolproof}
		By Theorem~\ref{theo:Fubini}, condition \ref{item:E} is equivalent to
		\begin{enumerate}[label=\ep{\normalfont{}i$'$}]
			\item\label{item:bis} $\dim (A \cap (X \times Y)) = \dim A$,
		\end{enumerate}
		and it follows from the \hyperref[theo:fib_dim]{dimension of fibers theorem} that \ref{item:bis} is equivalent to $\dim Y = \dim X$.
	\end{coolproof}
	
	\section{Expansion in algebraic bipartite graphs}\label{sec:proj}
	
	\subsection{Overview}\label{subsec:proj_intro}
	
	The central results of this section are Theorem~\ref{theo:proj} and its Corollary~\ref{corl:proj_sharp}. They play a key role in the proof of Theorem~\ref{theo:main} and are also interesting in their own right as natural extensions of the \hyperref[theo:fib_dim]{dimension of fibers theorem}.
	
	Consider an irreducible definable set $A \subseteq \F^n \times \F^m$ such that $\dim \proj_1 A = n$. Given an irreducible definable subset $X \subseteq \F^n$, what is the dimension of the fiber of $A$ over $X$, i.e., of the set $A \cap (X \times \F^m)$? For simplicity, assume that every component $B$ of $A \cap (X \times \F^m)$ is \emph{dominant}, i.e., $\dim \proj_1 B = \dim X$ \ep{the dimension of nondominant components is harder to control}. The ``expected'' answer is
	\begin{equation}\label{eq:DFT}
		\dim(A \cap (X \times \F^m)) \,=\, \dim A - n + \dim X,
	\end{equation}
	and we hope that \eqref{eq:DFT} holds for ``typical'' $X$. This hope is justified by the \hyperref[theo:fib_dim]{dimension of fibers theorem}. Indeed, let $C$ be the closure of the set of all $x \in \F^n$ such that $\dim A_x \neq \dim A - n$. By Theorem~\ref{theo:fib_dim}, $\dim C < n$, and if \eqref{eq:DFT} fails, then $X \subseteq C$. In other words, the \hyperref[theo:fib_dim]{dimension of fibers theorem} gives us a ``small'' definable set $C$ that contains every counterexample to \eqref{eq:DFT}.
	
\begin{figure}[h]
		\centering
		\begin{tikzpicture}
			\draw[dashed] (-0.7-1.2,0) -- node (M) [midway] {} (-0.3-2.5,2);
			
			\foreach \x in {1,...,10}
				\draw[dashed] (-0.7-1.2+\x*0.24,0) -- (-0.3-2.5+\x*0.5,2);
		
			\draw (0,0) circle [x radius=3, y radius=0.5];
			\draw (0,2) circle [x radius=4, y radius=0.7];
			
			\definecolor{light-gray}{gray}{0.95}
			
			\filldraw[fill=light-gray] (-0.7,0) circle [x radius=1.2, y radius=0.3];
			\filldraw[fill=light-gray] (-0.3,2) circle [x radius=2.5, y radius=0.5];
			
			\node at (-0.7,0) {$X$};
			\node at (-0.3,2) {$\proj_2 (A \cap (X \times \F^m))$};
			
			\node at (3.5,0) {$\F^n$};
			\node at (4.5,2) {$\F^m$};
			
			\node (L) at (-5,1) {$A \cap (X \times \F^m)$};
			
			\draw[-{Stealth[length=1.6mm]}] (L) to (M);
		\end{tikzpicture}
		\caption{A combinatorial interpretation of $\proj_2 (A \cap (X \times \F^m))$.}\label{fig:cartoon}
	\end{figure}
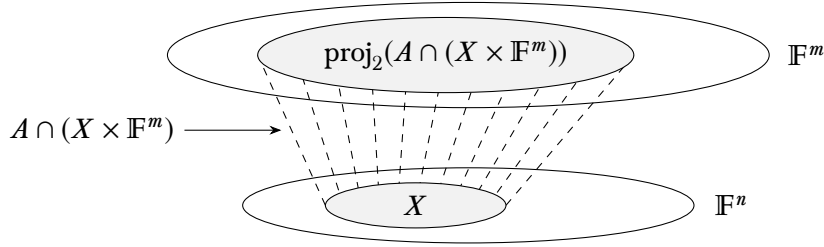

	We wish to obtain a version of this result for the dimension of the \emph{second projection} $\proj_2 (A \cap (X \times \F^m))$. Notice that $\proj_2 (A \cap (X \times \F^m))$ has a natural combinatorial interpretation: If we think of $A$ as the edge set of a bipartite graph with bipartition $(\F^n, \F^m)$, then $\proj_2 (A \cap (X \times \F^m))$ is the \emph{neighborhood} of $X$ in this graph \ep{see Fig.~\ref{fig:cartoon}}. There are two obvious upper bounds on $\dim \proj_2 (A \cap (X \times \F^m))$: First,
	\[
	\dim \proj_2 (A \cap (X \times \F^m)) \,\leq\, \dim (A \cap (X \times \F^m)) \,\stackrel{\mathclap{\eqref{eq:DFT}}}{=}\, \dim A - n + \dim X;
	\]
	and second, since $A \cap (X \times \F^m) \subseteq A$,
	\[
	\dim \proj_2 (A \cap (X \times \F^m)) \,\leq\, \dim \proj_2 A.
	\]
	The goal of this section is to show that, for ``typical'' $X$, at least one of these bounds must be tight:
	\begin{equation}\label{eq:dim_proj_intro}
		\dim \proj_2 (A \cap (X \times \F^m))  \,=\, \min \set{\dim A - n + \dim X, \, \dim \proj_2 A}.
	\end{equation}
	In other words, algebraic graphs are ``maximally expanding'': the dimension of the neighborhood of a ``typical'' set in such a graph is as large as it can possibly be.
	
	What is meant by a ``typical'' set here is a somewhat subtle issue. Based on the preceding discussion, it is tempting to conjecture that there should be a ``small'' definable set $C$ that contains every counterexample $X$ to \eqref{eq:dim_proj_intro}. This, however, need not be the case, as the following construction illustrates. Suppose that $n = 2$, $m = 1$, and $A \subset \F^2 \times \F$ is given by
	\begin{equation}\label{eq:example_lines}
		A \defeq \set{(x,y,z) \in \F^2 \times \F \,:\,  y = zx \text{ and } x \neq 0}.
	\end{equation}
	Then $A$ is an irreducible set of dimension $2$ and $\proj_2 A = \F$. Thus, for this $A$ and for any $1$-dimensional set $X \subset \F^2$, \eqref{eq:dim_proj_intro} turns into $\dim \proj_2 (A \cap (X \times \F))  = 1$. But if $X \subset \F^2$ is a straight line passing through the origin $(0,0)$, then $\dim \proj_2 (A \cap (X \times \F)) \leq 0$, and the union of all such lines is all of $\F^2$.
	
	Our approach is to allow the ``container'' $C$ to \emph{vary} with $X$, but in a very \emph{limited way}. \ep{This idea is inspired by the hypergraph containers method, where one builds a container for a given independent set $I$ based on a small fingerprint $F \subseteq I$.} To be more precise, let $\calF$ be a family of $k$-dimensional counterexamples to \eqref{eq:dim_proj_intro}. \ep{As in the statement of Theorem~\ref{theo:main}, we shall only work with \emph{parameterized} sets, so $\calF$ is really a family of \emph{definable functions} rather than sets, but for the purposes of the current informal discussion, this technicality may be ignored.} Imagine that Alice and Bob are playing the following game: Alice secretly chooses a set $X \in \calF$. Then she picks an $r$-dimensional subset $Y \subset X$, where $r < k$, and shows it to Bob. Bob's goal is to find, based on $Y$ alone, a ``container'' $C$ such that $\dim C < n$ and $X \subseteq C$. Corollary~\ref{corl:proj_sharp} asserts, roughly speaking, that there is a definable construction that allows Bob to win for all $X \in \calF$ and for generic $Y \subset X$. 
	
	As an illustration, consider again the set $A \subset \F^2 \times \F$ given by \eqref{eq:example_lines} and let $\calF$ be the family of all straight lines in $\F^2$ passing through the origin. Whatever $X \in \calF$ Alice chooses, when she shows Bob a generic point $y \in X$, he can win simply by making $C$ be the unique straight line passing through $(0,0)$ and $y$, in agreement with Corollary~\ref{corl:proj_sharp}.
	
	Before we proceed to formal statements and proofs, there is one more detail that is worth pointing out, namely what we mean by a \emph{generic} $r$-dimensional subset $Y \subset X$. As mentioned before, we only work with parametrized sets, so let $f \colon \F^k \rightharpoonup \F^n$ be the parameterization of $X$. We can then pick a generic affine map $\ell \in \L(r, k)$ and take $Y$ to be the subset of $X$ parameterized by the \emph{composition} $f \circ \ell \colon \F^r \rightharpoonup \F^n$. The reader may be wondering why we restrict our attention to \emph{affine} maps $\ell$ instead of allowing rational maps of arbitrary degree. The reason, roughly speaking, is that we wish the subset $Y \subset X$ to contain ``less information'' than $X$ itself. When $\ell$ is affine, this intuition is justified as $f \circ \ell$ is a map of the {same degree} as $f$ but in fewer variables. Note that this construction can be iterated; in other words, given $\lambda \in \L(s, r)$, we can look at the triple composition $f \circ \ell \circ \lambda \colon \F^s \rightharpoonup \F^n$, which can be interpreted as picking a generic $s$-dimensional subset of a generic $r$-dimensional subset of $X$. Unsurprisingly, this operation is essentially equivalent to directly picking a generic $s$-dimensional subset of $X$ \ep{see Lemma~\ref{lemm:comp}; this again relies on the fact that the maps $\ell$ and $\lambda$ are affine, so their composition $\ell \circ \lambda$ is affine as well}, and this simple fact will be crucial for our arguments \ep{in particular, in the proof of Proposition~\ref{prop:reduce}}.
	

	\subsection{Definable families of functions}
	
	By a \emphd{definable family of functions} from $\F^k$ to $\F^n$ we mean a pair $(\calF, \eval)$,
	where $\calF$ is a  definable set (in some power of $\F$) and $\eval \colon \calF \times \F^k \rightharpoonup \F^n$ is a definable partial function, called the \emphd{evaluation map}, such that $\forall^\ast f \in \calF \, \forall^\ast y \in \F^k\, (\eval(f, y) \text{ is defined})$.
	The evaluation map is usually clear from the context, so we omit it and simply write $\calF$ instead of $(\calF, \eval)$. We write $\calF \colon \F^k \Rrightarrow \F^n$ to indicate that $\calF$ is a definable family of functions from $\F^k$ to $\F^n$. Each $f \in \calF$ gives rise to the partial map
	\[
		\F^k \rightharpoonup \F^n \colon y \mapsto \eval(f, y),
	\]
	which we also denote by $f$ (so expressions like $f(y)$ or $f \circ g$ must be interpreted accordingly). The basic example of a definable family of functions from $\F^k$ to $\F^n$ is $\Rat_d(k,n;q)$ for a nonzero polynomial $q \in \F[x_1, \ldots, x_k]$ (with the natural evaluation map).
	
	\begin{defn}\label{defn:comprehensive}	
	A definable family of functions $\calF \colon \F^k \Rrightarrow \F^n$ is \emphd{comprehensive} if the following statements are equivalent for every definable set $C \subseteq \F^n$:
	\begin{enumerate}[label=\ep{\normalfont{}C\arabic*}]
		\item\label{item:C1} $\forall^\ast f \in \calF \quad \forall^\ast y \in \F^k \quad (f(y) \in C)$; 
		
		\item\label{item:C2} $\forall^\ast x \in \F^n\, (x \in C)$, i.e., $\dim C = n$. 
	\end{enumerate}
	\end{defn}
	
	For instance, the family $\Rat_d(k,n;q)$ is comprehensive, since for all $y \in \F^k$ with $q(y) \neq 0$ and for all $x$, $x' \in \F^n$, $\dim \set{f \in \Rat_d(k,n;q) \,:\, f(y) = x} = \dim \set{f \in \Rat_d(k,n;q) \,:\, f(y) = x'}$; indeed, the map \[\set{f \in \Rat_d(k,n;q) \,:\, f(y) = x} \to \set{f \in \Rat_d(k,n;q) \,:\, f(y) = x'} \colon f \mapsto f + q(y)(x'-x)/q\] is a definable bijection.
	
	\begin{prop}\label{prop:C1toC2}
		Let $\calF \colon \F^k \Rrightarrow \F^n$ be a definable family of functions. If $\calF$ is irreducible and satisfies the implication \ref{item:C1} $\Longrightarrow$ \ref{item:C2}, then $\calF$ is comprehensive.
	\end{prop}
	\begin{coolproof}
		We have to prove \ref{item:C2} $\Longrightarrow$ \ref{item:C1}. To that end, let $C \subseteq \F^n$ be a definable set for which \ref{item:C1} fails. Since $\calF$ is irreducible, the negation of \ref{item:C1} is equivalent to $\forall^\ast f \in \calF \, \forall^\ast y \in \F^k \, (f(y) \not \in C)$, so applying \ref{item:C1} $\Longrightarrow$ \ref{item:C2} with $\F^n \setminus C$ in place of $C$ yields $\dim (\F^n \setminus C) = n$, i.e., $\dim C < n$, as desired.
	\end{coolproof}
	
	The following construction will be useful. Let $\calF \colon \F^k \Rrightarrow \F^n$ and $\calG \colon \F^r \Rrightarrow \F^k$ and assume that $\calG$ is comprehensive. Define $\calF \otimes \calG \colon \F^r \Rrightarrow \F^n$, the \emphd{composition} of $\calF$ and $\calG$, as follows: As a set, $\calF \otimes \calG$ is equal to $\calF \times \calG$, and the evaluation map on $\calF\otimes \calG$ is given by
	\[
	(f,g)(z) \defeq (f \circ g)(z) \text{ for all } f \in \calF,\ g \in \calG, \text{ and } z \in \F^r.
	\]
	
	\begin{prop}
		The above definition is correct; that is, for all $\calF \colon \F^k \Rrightarrow \F^n$ and $\calG \colon \F^r \Rrightarrow \F^k$, if $\calG$ is comprehensive, then $\calF \otimes \calG$ is a definable family of functions from $\F^r$ to $\F^n$.
	\end{prop}
	\begin{coolproof}
		Due to Corollary~\ref{corl:FFubini}, we have to verify that
		\[
		\forall^\ast f \in \calF \quad \forall^\ast g \in \calG \quad \forall^\ast z \in \F^r \qquad (f \circ g)(z) \text{ is defined}.
		\]
		Since the family $\calG$ is comprehensive, this is equivalent to $\forall^\ast f \in \calF \, \forall^\ast y \in \F^k \, (f(y) \text{ is defined})$, which holds by definition, since $\calF$ is a definable family of functions.
	\end{coolproof}
	
	\begin{prop}\label{prop:compr}
		Let $\calF \colon \F^k \Rrightarrow \F^n$ and $\calG \colon \F^r \Rrightarrow \F^k$. If $\calF$ and $\calG$ are comprehensive, then so is $\calF \otimes \calG$.
	\end{prop}
	\begin{coolproof}
		For any definable set $C \subseteq \F^n$, we have
		\begin{align*}
		&\forall^\ast f \in \calF \quad \forall^\ast g \in \calG \quad \forall^\ast z \in \F^r \quad (f \circ g)(z) \in C \\
		[\text{$\calG$ is comprehensive}] \quad \Longleftrightarrow \quad &\forall^\ast f \in \calF \quad \forall^\ast y \in \F^k \quad f(y) \in C \\
		[\text{$\calF$ is comprehensive}] \quad  \Longleftrightarrow \quad &\forall^\ast x \in \F^n \quad x \in C. \qedhere
		\end{align*}
	\end{coolproof}
	
	\subsection{Containers}
	
	A crucial concept for the results of this section is that of an \emph{$r$-container}. Informally, an $r$-container $\mathfrak{C}$ for a definable family $\calF \colon \F^k \Rrightarrow \F^n$ is a definable rule that, given $f \in \calF$ and $\ell \in \L(r, k)$, outputs a subset $\mathfrak{C}_{f \circ \ell} \subseteq \F^n$ that only depends on the {composition} $f \circ \ell$ and not on $f$ and $\ell$ themselves. Here is the precise definition:
	
	\begin{defn}
		Let $\calF \colon \F^k \Rrightarrow \F^n$ and let $r <k$. An \emphd{$r$-container} for $\calF$ is a definable set
	\[
		\mathfrak{C} \,\subseteq\, \calF \times \L(r, k) \times \F^n
	\]
	such that for all $f$, $g \in \calF$ and $\ell$, $\lambda \in \L(r,k)$ with $f \circ \ell = g \circ \lambda$, and for all $x \in \F^n$, we have
	\[
	(f, \ell, x) \in \mathfrak{C} \quad \Longleftrightarrow \quad (g, \lambda, x) \in \mathfrak{C}.
	\]
	If $\mathfrak{C}$ is an $r$-container for $\calF$, then, for any definable map $\phi \colon \F^r \rightharpoonup \F^n$, we write
	\[
	\mathfrak{C}_\phi \,\defeq\, \set{x \in \F^n \,:\, (f, \ell, x) \in \mathfrak{C} \text{ for some } f \in \calF \text{ and } \ell \in \L(r,k) \text{ such that } f \circ \ell = \phi}.
	\]
	\end{defn}
	
	
	\begin{defn}\label{defn:NCP}
		Let $\calF \colon \F^k \Rrightarrow \F^n$ and $r \in \set{-1,\, 0,\, 1,\, \ldots,\, k-1}$. The family $\calF$ is $r$-\emphd{uncontainable} if either $r = - 1$ and $\calF$ is comprehensive, or else, $r \in \N$ and the following statements are equivalent for every $r$-container $\mathfrak{C}$ for $\calF$:
	\begin{enumerate}[label=\ep{\normalfont{}U\arabic*}]
		\item\label{item:U_capture} $\forall^\ast f \in \calF \quad \forall^\ast \ell \in \L(r, k) \quad \forall^\ast y \in \F^k \quad (f(y) \in \mathfrak{C}_{f \circ \ell})$;
		
		\item\label{item:U_small} $\forall^\ast f \in \calF \quad \forall^\ast \ell \in \L(r, k) \quad \forall^\ast x \in \F^n \quad (x \in \mathfrak{C}_{f \circ \ell})$.
	\end{enumerate}
	\end{defn}
	
	The implication \ref{item:U_capture} $\Longrightarrow$ \ref{item:U_small} in Definition~\ref{defn:NCP} can be informally summarized as, ``It is impossible to put most elements of $\calF$ into small containers,'' hence the term ``uncontainable.'' The opposite implication \ref{item:U_small} $\Longrightarrow$ \ref{item:U_capture} says that, conversely, ``large containers must capture most elements of $\calF$.'' It is convenient to include the second implication as part of Definition~\ref{defn:NCP}, even though in most cases of interest to us it will follow automatically, see Proposition~\ref{prop:U1toU2}.
	
	When working with $r$-uncontainable families, we often rely on the following basic fact:
	
	\begin{lemma}\label{lemm:comp}
		Suppose that $s \leq k$, $s \leq r$ and let $\lambda \in \L(s, r)$ be an injective affine map. Then, for every definable set $A \subseteq \L(s, k)$, we have
		\[
			\forall^\ast \ell \in \L(r, k) \quad \ell \circ \lambda \in A \qquad \Longleftrightarrow \qquad \forall^\ast \phi \in \L(s,k) \quad \phi \in A.
		\]
	\end{lemma}
	\begin{coolproof}
		Since $s \leq k$, a generic affine map $\phi \in \L(s, k)$ is injective, and, similarly, $\ell \circ \lambda$ is injective for generic $\ell \in \L(r, k)$. It remains to observe that for all injective $\phi \in \L(s, k)$, the dimension of the set $\set{\ell \in \L(r, k) \,:\, \ell \circ \lambda = \phi}$ is the same.
	\end{coolproof}
	
	\begin{prop}
		Let $\calF \colon \F^k \Rrightarrow \F^n$ and $s$, $r \in \set{-1,\, 0,\, 1,\, \ldots,\, k-1}$. If $s \leq r$, then
	\[
		\text{$\calF$ is $r$-uncontainable} \;\Longrightarrow\; \text{$\calF$ is $s$-uncontainable}.
	\]
	\end{prop}
	\begin{coolproof}
		Suppose that $0 \leq s \leq r < k$ and $\calF$ is $r$-uncontainable \ep{the proof in the case $s = -1$ is similar, and we omit it}. Let $\mathfrak{C}$ be an arbitrary $s$-container for $\calF$. Fix any injective map $\lambda \in \L(s, r)$ and define an $r$-container $\tilde{\mathfrak{C}}$ by $\tilde{\mathfrak{C}}_{f \circ \ell} \defeq \mathfrak{C}_{f \circ \ell \circ \lambda}$ for all $f \in \calF$ and $\ell \in \L(r, k)$.
		We claim that \ref{item:U_capture} and \ref{item:U_small} for $\mathfrak{C}$ are equivalent to the corresponding statements for $\tilde{\mathfrak{C}}$, which implies that $\calF$ is $s$-uncontainable. Indeed, \ref{item:U_capture} for $\tilde{\mathfrak{C}}$ takes the form
		\[
			\forall^\ast f \in \calF \quad \forall^\ast \ell \in \L(r, k) \quad \forall^\ast y \in \F^k \qquad f(y) \in \mathfrak{C}_{f \circ \ell \circ \lambda}.
		\]
		By Lemma~\ref{lemm:comp}, this is equivalent to
		\[
			\forall^\ast f \in \calF \quad \forall^\ast \phi \in \L(s, k) \quad \forall^\ast y \in \F^k \qquad f(y) \in \mathfrak{C}_{f \circ \phi},
		\]
		which is precisely \ref{item:U_capture} for $\mathfrak{C}$. The argument for \ref{item:U_small} is similar.
	\end{coolproof}
	
	\begin{prop}\label{prop:U1toU2}
		Let $\calF \colon \F^k \Rrightarrow \F^n$ and $r \in \set{-1,\, 0,\, 1,\, \ldots,\, k-1}$. If $\calF$ is irreducible and satisfies the implication \ref{item:U_capture} $\Longrightarrow$ \ref{item:U_small}, then $\calF$ is $r$-uncontainable.
	\end{prop}
	\begin{coolproof}
		The proof is the same, \emph{mutatis mutandis}, as the proof of Proposition~\ref{prop:C1toC2}.
	\end{coolproof}
	
	The family $\L(k,n)$ is easily seen to be $(k-1)$-uncontainable \ep{this is a special case of the more general Proposition~\ref{prop:no_cont} proved in Section~\ref{sec:proof}}. This fact yields the following extension of Proposition~\ref{prop:compr}:
	
	\begin{prop}\label{prop:reduce}
		Let $\calF \colon \F^k \Rrightarrow \F^n$ and $r \in \set{-1,\, 0,\, 1,\, \ldots,\, k-1}$. Suppose that $s \leq k$ and let $t \defeq \min \set{r, s-1}$. If $\calF$ is $r$-uncontainable, then $\calF \otimes \L(s, k)$ is $t$-uncontainable.
	\end{prop}
	\begin{coolproof}
		The case $t = -1$ is handled by Proposition \ref{prop:compr}, so, from now on, assume that $t \in \N$.
		
		Let $\mathfrak{C}$ be a $t$-container for $\calF \otimes \L(s, k)$. Then \ref{item:U_capture} takes the form
		\[
			\forall^\ast f \in \calF \quad \forall^\ast \ell \in \L(s,k) \quad \forall^\ast \lambda \in \L(t, s) \quad \forall^\ast z \in \F^s \qquad (f \circ \ell)(z) \in \mathfrak{C}_{f \circ \ell \circ \lambda}.
		\]
		Since $\L(s, k)$ is $(s-1)$-uncontainable and $t \leq s-1$, this yields
		\[
		\forall^\ast f \in \calF \quad \forall^\ast \ell \in \L(s,k) \quad \forall^\ast \lambda \in \L(t, s) \quad \forall^\ast y \in \F^k \qquad f(y) \in \mathfrak{C}_{f \circ \ell \circ \lambda},
		\]
		which, by Lemma~\ref{lemm:comp}, is equivalent to
		\[
		\forall^\ast f \in \calF \quad \forall^\ast \phi \in \L(t, k) \quad \forall^\ast y \in \F^k \qquad f(y) \in \mathfrak{C}_{f \circ \phi}.
		\]
		Since $\calF$ is $r$-uncontainable and $t \leq r$, we conclude that
		\[
		\forall^\ast f \in \calF \quad \forall^\ast \phi \in \L(t, k) \quad \forall^\ast x \in \F^n \qquad x \in \mathfrak{C}_{f \circ \phi},
		\]
		which, by Lemma~\ref{lemm:comp} again, can be rewritten as
		\[
			\forall^\ast f \in \calF \quad \forall^\ast \ell \in \L(s,k) \quad \forall^\ast \lambda \in \L(t, s) \quad \forall^\ast x \in \F^n \qquad x \in \mathfrak{C}_{f \circ \ell \circ \lambda}.
		\]
		But the last expression is exactly the relevant instance of \ref{item:U_small}.
	\end{coolproof}

	\subsection{Typical dimension of neighborhoods}\label{subsec:proj}
	
	Let $A \subseteq \F^n \times \F^m$ be a definable set. We say that $A$ is \emphd{dominant} if $\dim\proj_1 A = n$. For a definable map $\phi \colon \F^k \rightharpoonup \F^n$, the \emphd{fiber of $A$ over $\phi$} is the set $A_\phi \subseteq \F^k \times \F^m$ given by
	\[
		(y_1,\, x_2) \in A_\phi \quad \vcentcolon\Longleftrightarrow \quad (\phi(y_1),\, x_2) \in A.
	\]
	The following lemma is a precise statement of the corollary of the \hyperref[theo:fib_dim]{dimension of fibers theorem} discussed in the beginning of \S\ref{subsec:proj_intro}:
	
	\begin{lemma}\label{lemma:dominant}
		Let $A \subseteq \F^n \times \F^m$ be a dominant definable set and let $\calF \colon \F^k \Rrightarrow \F^n$. If $\calF$ is comprehensive, then, for generic $f \in \calF$, the set $A_f \subseteq \F^k \times \F^m$ is dominant and the dimension of every dominant irreducible component of $A_f$ is $\dim A - n + k$.
	\end{lemma}
	
	Note that, by Lemma~\ref{lemma:comps_are_definable}, the set of all $f \in \calF$ that satisfy the conclusion of Lemma~\ref{lemma:dominant} is definable.
	
	\begin{coolproof}
		Let $C \subseteq \F^n$ be the set of all $x \in \F^n$ such that $A_x \neq \0$ and the dimension of every irreducible component of $A_x$ is $\dim A - n$. Since $A$ is dominant, by Theorem~\ref{theo:fib_dim}, $\dim C = n$. Applying the implication \ref{item:C2} $\Longrightarrow$ \ref{item:C1}, we obtain that for generic $f \in \calF$ and for generic $y \in \F^k$, $f(y) \in C$; in other words, $(A_f)_y \neq \0$ and the dimension of every irreducible component of $(A_f)_y$ is $\dim A - n$. Since $(A_f)_y \neq \0$ for generic $y \in \F^k$, we conclude that $A_f$ is dominant. Let $B \subseteq A_f$ be any dominant component of $A_f$. For generic $y \in \F^k$, we have $\dim B_y = \dim B - k$. Since $B_y \subseteq (A_f)_y$, this implies that $\dim B \leq \dim A - n + k$. Suppose that $\dim B < \dim A - n + k$. Then, for generic $y \in \F^k$, $\dim B_y$ is less than the dimension of every irreducible component of $(A_f)_y$, and hence $B_y$ is contained in the closure of $(A_f)_y \setminus B_y$. This means that $B \subseteq \cl{A_f \setminus B}$, contradicting the fact that $B$ is a component of $A_f$.
	\end{coolproof}
	
	We are now ready to state and prove the central result of this section:
	
	\begin{theo}\label{theo:proj}
		Let $A \subseteq \F^n \times \F^m$ be a dominant irreducible definable set and let $\calF \colon \F^k \Rrightarrow \F^n$ be an irreducible definable family of functions. Let
		\[
			r \,\defeq\, \min \set{k, \, \dim\proj_2 A + n - \dim A} \,-\, 1.
		\]
		If $\calF$ is $r$-uncontainable, then, for generic $f \in \calF$, every dominant component $B$ of $A_f$ satisfies
		\[
			\dim \proj_2 B \,=\, \min\set{\dim A - n + k,\, \dim \proj_2 A}.
		\]
	\end{theo}
	
	
	\begin{coolproof}
		The proof is by induction on $k$. The base case $k = 0$ is a restatement of Theorem~\ref{theo:fib_dim}, so suppose that $k \geq 1$. Let $\calF \colon \F^k \Rrightarrow \F^n$ be an $r$-uncontainable irreducible definable family. It follows from Lemma~\ref{lemma:dominant} that, for generic $f \in \calF$, every dominant component $B$ of $A_f$ satisfies
		\begin{equation}\label{eq:upper_bound}
			\dim \proj_2 B \,\leq\, \min\set{\dim B,\, \dim \proj_2 A} \,=\, \min\set{\dim A - n + k,\, \dim \proj_2 A}.
		\end{equation}
		Suppose that the conclusion of Theorem~\ref{theo:proj} fails for $\calF$. Due to \eqref{eq:upper_bound} and since $\calF$ is irreducible, this means that for generic $f \in \calF$, the set $A_f$ has a dominant component $B$ with
		\begin{equation}\label{eq:less}
			\dim \proj_2 B \,<\, \min\set{\dim A - n + k,\, \dim \proj_2 A}.
		\end{equation}
		Consider the family $\calF \otimes \L(k-1, k) \colon \F^{k-1} \Rrightarrow \F^n$. Since the families $\calF$ and $\L(k-1, k)$ are comprehensive, so is $\calF \otimes \L(k-1, k)$ by Proposition~\ref{prop:compr}. By Lemma~\ref{lemma:dominant}, this implies that for generic $f \in \calF$ and $\ell \in \L(k-1, k)$, every dominant component $C$ of $A_{f \circ \ell}$ satisfies
		\begin{equation}\label{eq:C}
			\dim C \,=\, \dim A - n + k-1.
		\end{equation}
		Furthermore, we can use the inductive assumption and apply Theorem~\ref{theo:proj} to $\calF \otimes \L(k-1, k)$. Set \[s \,\defeq\, \min \set{k-1, \, \dim\proj_2 A + n - \dim A} \,-\, 1.\] Then $s = \min \set{r, k-2}$, and, since $\calF$ is $r$-uncontainable, Proposition~\ref{prop:reduce} yields that $\calF \otimes \L(k-1, k)$ is $s$-uncontainable. Hence, the conclusion of Theorem~\ref{theo:proj} holds for $\calF \otimes \L(k-1, k)$; in other words, for generic $f \in \calF$ and $\ell \in \L(k-1, k)$, every dominant component $C$ of $A_{f \circ \ell}$ satisfies
		\begin{equation}\label{eq:Cproj}
			\dim \proj_2 C \,=\, \min\set{\dim A - n + k - 1,\, \dim \proj_2 A}.
		\end{equation}
		
		\stepcounter{ForClaims} \renewcommand{\theForClaims}{\ref{theo:proj}}
		\begin{claim}\label{claim:descent}
			The following statements are valid:
			\begin{enumerate}[label=\ep{\itshape{\alph*}}, wide]
				\item\label{item:descent_a} $\dim \proj_2 A \geq \dim A - n + k$, and hence $r=k-1$;
				\item\label{item:descent_b} for generic $f \in \calF$ and $\ell \in \L(k-1, k)$, there exist a dominant component $B$ of $A_f$ and a dominant component $C$ of $A_{f \circ \ell}$ such that $\cl{\proj_2 B} = \cl{\proj_2 C}$ and $\dim \proj_2 B = \dim \proj_2 C = \dim A - n + k - 1$. 
			\end{enumerate}
		\end{claim}
		\begin{claimproof}
		Consider generic $f \in \calF$ and $\ell \in \L(k-1, k)$. Let $B$ be an arbitrary dominant component of $A_f$ satisfying \eqref{eq:less}. Recall that, by Lemma~\ref{lemma:dominant}, $\dim B = \dim A - n + k$. Since the family $\L(k-1, k)$ is comprehensive, we may apply Lemma~\ref{lemma:dominant} to $B$ in order to conclude that the set $B_\ell$ is dominant and the dimension of every dominant component of $B_\ell$ is $\dim B - k + (k-1) = \dim A - n + k - 1$.
		Let $B'$ be any dominant component of $B_\ell$. Since $B_\ell \subseteq A_{f \circ \ell}$, there is a \ep{necessarily also dominant} component $C$ of $A_{f \circ \ell}$ such that $B' \subseteq C$. From \eqref{eq:C}, we see that $\dim B' = \dim C$, and hence the closures of $B'$ and $C$ coincide. By Corollary~\ref{corl:proj_cl}, projections commute with closures, and therefore
		\[
		\cl{\proj_2 B'} = \cl{\proj_2 \cl{B'}} = \cl{\proj_2 \cl{C}} = \cl{\proj_2 C}, \qquad \text{so} \qquad \dim \proj_2 B' = \dim \proj_2 C.
		\]
		Combining this with \eqref{eq:less} and \eqref{eq:Cproj}, we obtain the following chain of \ep{in}equalities:
		\begin{align}
			\min\set{\dim A - n + k,\, \dim \proj_2 A} \,&>\, \dim \proj_2 B \,\geq\, \dim \proj_2 B_\ell \,\geq\, \dim \proj_2 B' \nonumber\\
			&=\, \dim \proj_2 C \,=\, \min\set{\dim A - n + k - 1,\, \dim \proj_2 A}. \label{eq:ineq}
		\end{align}
		By comparing the first and the last terms in \eqref{eq:ineq}, we see that $\dim \proj_2 A \geq \dim A - n + k$ and \[\dim \proj_2 B \,=\, \dim \proj_2 C \,=\, \dim A - n + k -1.\] Since $\proj_2 B$ is irreducible and $\cl{\proj_2 B} \supseteq \cl{\proj_2 C}$, this yields $\cl{\proj_2 B} = \cl{\proj_2 C}$, and we are done.
		\end{claimproof}
		
		Now we define a $(k-1)$-container $\mathfrak{C}$ for $\calF$ as follows: For each $f \in \calF$ and $\ell \in \L(k-1, k)$, let $\mathfrak{C}_{f \circ \ell}$ be the set of all $x \in \F^n$ such that the set $A_{f \circ \ell}$ has a dominant component $C$ with
		\[
			\dim \proj_2 C = \dim A - n + k -1 \qquad \text{and} \qquad \dim (A_x \cap \proj_2 C) \geq \dim A - n.
		\]
		Note that $\mathfrak{C}$ is definable due to Lemma~\ref{lemma:comps_are_definable} \ep{and other results of \S\ref{subsec:model_theory}}.
		
		\begin{claim}\label{claim:U1}
			The $(k-1)$-container $\mathfrak{C}$ satisfies \ref{item:U_capture}; in other words, we have
			\[
				\forall^\ast f \in \calF \quad \forall^\ast \ell \in \L(k-1,k) \quad \forall^\ast y \in \F^k \qquad f(y) \in \mathfrak{C}_{f \circ \ell}.
			\]
		\end{claim}
		\begin{claimproof}
			Pick generic $f \in \calF$ and $\ell \in \L(k-1, k)$ and let $B$ and $C$ be dominant components of $A_f$ and $A_{f \circ \ell}$ respectively with the properties specified by Claim~\ref{claim:descent}\ref{item:descent_b}. Since $\dim \proj_2 C \,=\, \dim A - n + k - 1$, every $y \in \F^k$ with $f(y) \not \in \mathfrak{C}_{f \circ \ell}$ must satisfy
			\begin{equation}\label{eq:small_intersection}
				\dim (B_y \cap \proj_2 C) \,\leq\, \dim ((A_f)_y \cap \proj_2 C) \,<\, \dim A - n.
			\end{equation}
			We claim that \eqref{eq:small_intersection} fails for generic $y \in \F^k$, which gives the desired conclusion. Indeed,
			we have $\dim B = \dim A - n + k$, so, for generic $y \in \F^k$, $\dim B_y = \dim A - n$, and if such $y$ satisfies \eqref{eq:small_intersection}, then
			\begin{equation}\label{eq:large_complement}
				\dim (B_y \setminus \proj_2 C) \,=\, \dim A - n.
			\end{equation}
			If \eqref{eq:large_complement} holds for generic $y \in \F^k$, then
			\[
				\dim (B \setminus (\F^k \times \proj_2 C)) \,\geq\, \dim A - n + k \,=\, \dim B.
			\]
			Since $B$ is irreducible, Corollary~\ref{corl:proj_cl} then yields
			\[
				\dim (\proj_2 B \setminus \proj_2 C) \,=\, \dim \proj_2 (B \setminus (\F^k \times \proj_2 C)) \,=\, \dim \proj_2 B,
			\]
			which contradicts the fact that $\cl{\proj_2 B} = \cl{\proj_2 C}$.
		\end{claimproof}
		
		Since, by Claim~\ref{claim:descent}\ref{item:descent_a}, the family $\calF$ is $(k-1)$-uncontainable, we deduce from Claim~\ref{claim:U1} that
		\[
			\forall^\ast f \in \calF \quad \forall^\ast \ell \in \L(k-1,k) \quad \forall^\ast x \in \F^n \qquad x \in \mathfrak{C}_{f \circ \ell}.
		\] 
		Take any $f \in \calF$ and $\ell \in \L(k-1, k)$ such that $\forall^\ast x\in \F^n \,(x \in \mathfrak{C}_{f \circ \ell})$. The set $A_{f \circ \ell}$ has only finitely many components, so we can choose a dominant component $C$ of $A_{f \circ \ell}$ so that
		\[
			\dim \proj_2 C = \dim A - n + k -1 \qquad \text{and} \qquad \forall^\ast x \in \F^n \, (\dim (A_x \cap \proj_2 C) \geq \dim A - n).
		\]
		Then we have $\dim (A \cap (\F^n \times \proj_2 C)) \geq \dim A$, and, since $A$ is irreducible, Corollary~\ref{corl:proj_cl} yields
		\[
			\dim \proj_2 C \,=\, \dim \proj_2 (A \cap (\F^n \times \proj_2 C)) \,=\, \dim \proj_2 A.
		\]
		But by Claim~\ref{claim:descent}\ref{item:descent_a}, $\dim \proj_2 A \geq \dim A -n + k > \dim \proj_2 C$; a contradiction.
	\end{coolproof}
	
	To establish a connection between Theorem~\ref{theo:proj} and the motivating discussion in \S\ref{subsec:proj_intro}, consider the contrapositive of Theorem~\ref{theo:proj}. Let $A \subseteq \F^n \times \F^m$ be a dominant irreducible definable set and suppose that $\calF \colon \F^k \Rrightarrow \F^n$ is an irreducible definable family consisting of functions $f$ that \emph{violate} the conclusion of Theorem~\ref{theo:proj}. Such a family $\calF$ cannot be $r$-uncontainable, and hence, by Proposition~\ref{prop:U1toU2}, there exists an $r$-container $\mathfrak{C}$ for $\calF$ such that for generic $f \in \calF$ and $\ell \in \L(r,k)$, we have
	\[\dim \mathfrak{C}_{f \circ \ell} < n, \qquad \text{yet} \qquad \dim f^{-1}(\mathfrak{C}_{f \circ \ell}) = k.\] In other words, based on $f \circ \ell$ alone, we can definably build a ``small'' set $\mathfrak{C}_{f \circ \ell}$ that contains ``most'' of $f(\F^k)$. This result can be further strengthened in two ways: first, we can actually do this for \emph{all} $f \in \calF$ \ep{and not just for a \emph{generic} subset}; second, the family $\calF$ need not be irreducible:
	
	\begin{corl}\label{corl:proj_sharp}
		Let $A \subseteq \F^n \times \F^m$ be a dominant irreducible definable set and let \[r \,\defeq\, \min \set{k,\, \dim\proj_2 A + n - \dim A} \,-\, 1.\]
		Let $\calF \colon \F^k \Rrightarrow \F^n$ be a definable family of functions such that for all $f \in \calF$, $\dim \dom(f) = k$. Suppose that for every $f \in \calF$, either the set $A_f$ is not dominant, or else, $A_f$ has a dominant component $B$ with
		\[
			\dim \proj_2 B \,\neq\, \min\set{\dim A - n + k,\, \dim \proj_2 A}.
		\]
		Then there is an $r$-container $\mathfrak{C}$ for $\calF$ such that for all $f \in \calF$ and for generic $\ell \in \L(r, k)$, we have
		\[\dim \mathfrak{C}_{f \circ \ell} < n, \qquad \text{yet} \qquad \dim f^{-1}(\mathfrak{C}_{f \circ \ell}) = k.\]
	\end{corl}
	\begin{coolproof}
		For fixed $A \subseteq \F^n \times \F^m$ and $k$, let $\calF \colon \F^k \Rrightarrow \F^n$ be a counterexample to Corollary~\ref{corl:proj_sharp} that minimizes $\dim \calF$ and has the fewest irreducible components among all counterexamples with dimension equal to $\dim \calF$. Clearly, $\calF \neq \0$. Since for all $f \in \calF$, $\dim \dom(f) = k$, every definable subset of $\calF$ forms a definable family of functions with the evaluation map inherited from $\calF$. Call a definable subfamily $\calF' \subseteq \calF$ \emphd{small} if there is an $r$-container $\mathfrak{C}$ for $\calF'$ such that for all $f \in \calF'$ and for generic $\ell \in \L(r, k)$, we have $\dim \mathfrak{C}_{f \circ \ell} < n$ but $\dim f^{-1}(\mathfrak{C}_{f \circ \ell}) = k$. By assumption, $\calF$ itself is not small.
		
		\stepcounter{ForClaims} \renewcommand{\theForClaims}{\ref{corl:proj_sharp}}
		\begin{claim}\label{claim:union}
			If $\calF_1$, $\calF_2 \subseteq \calF$ are small definable subfamilies, then so is $\calF_1 \cup \calF_2$.
		\end{claim}
		\begin{claimproof}
			For each $i \in [2]$, let $\mathfrak{C}_i$ be an $r$-container for $\calF_i$ such that for all $f \in \calF_i$ and for generic $\ell \in \L(r, k)$, we have $\dim (\mathfrak{C}_i)_{f \circ \ell} < n$ and $\dim f^{-1}((\mathfrak{C}_i)_{f \circ \ell}) = k$. We may in fact assume that $\dim (\mathfrak{C}_i)_{f \circ \ell} < n$ for \emph{all} $\ell \in \L(r, k)$, since otherwise we can replace $\mathfrak{C}_i$ with $\mathfrak{C}_i'$ given by
			\[
			(\mathfrak{C}_i')_{f \circ \ell} \defeq \begin{cases}
			(\mathfrak{C}_i)_{f \circ \ell} &\text{if } \dim (\mathfrak{C}_i)_{f \circ \ell} < n;\\
			\0 &\text{otherwise}.
			\end{cases}
			\]
			Define an $r$-container $\mathfrak{C}$ for $\calF_1 \cup \calF_2$ by setting $\mathfrak{C}_{f \circ \ell} \defeq (\mathfrak{C}_1)_{f \circ \ell} \cup (\mathfrak{C}_2)_{f \circ \ell}$. It is easy to see that $\mathfrak{C}$ has all the desired properties.
		\end{claimproof}
		
		We claim that $\calF$ is irreducible. Otherwise, we can write $\calF = \calF_1 \cup \calF_2$, where each of $\calF_1$, $\calF_2$ is a definable family of functions that has fewer irreducible components than $\calF$. By the choice of $\calF$, the families $\calF_1$ and $\calF_2$ are small, and hence so is $\calF$ by Claim~\ref{claim:union}, which is a contradiction.
		
		If $\calF$ were $r$-uncontainable, then we would be able to use Lemma~\ref{lemma:dominant} and Theorem~\ref{theo:proj} to conclude that for generic $f \in \calF$, the set $A_f$ is dominant and every dominant component $B$ of $A_f$ satisfies
		\[
		\dim \proj_2 B \,=\, \min\set{\dim A - n + k,\, \dim \proj_2 A}.
		\]
		This conclusion contradicts the assumptions on $\calF$, so $\calF$ cannot be $r$-uncontainable. From Proposition~\ref{prop:U1toU2}, it follows that there exists an $r$-container $\mathfrak{C}$ for $\calF$ that satisfies \ref{item:U_capture} but fails \ref{item:U_small}; in other words, for generic $f \in \calF$ and $\ell \in \L(r, k)$, we have $\dim f^{-1}(\mathfrak{C}_{f \circ \ell}) = k$ but $\dim \mathfrak{C}_{f \circ \ell} < n$. Let
		\[
			\calF_1 \defeq \set{f \in \calF \,:\, \forall^\ast \ell \in \L(r, k) \, (\dim f^{-1}(\mathfrak{C}_{f \circ \ell}) = k \text{ and } \dim \mathfrak{C}_{f \circ \ell} < n)} \qquad \text{and} \qquad \calF_2 \defeq \calF \setminus \calF_1.
		\]
		The $r$-container $\mathfrak{C}$ certifies that $\calF_1$ is small. But $\dim \calF_2 < \dim \calF$, so $\calF_2$ is also small by the choice of $\calF$. Claim~\ref{claim:union} then implies that $\calF$ itself is small---a contradiction.
	\end{coolproof}

	\section{Proof of Theorem~\ref{theo:main}}\label{sec:proof}
	
	\subsection{Almost density is necessary}\label{subsec:lower_bound}
	
	The following simple linear-algebraic fact will be rather useful \ep{it is partially responsible for the lower bound $d \geq t-1$ in the statement of Theorem~\ref{theo:main}}:
	
	\begin{lemma}\label{lemma:injective}
		Let $t \geq 1$. Fix $d \geq t-1$ and a nonzero polynomial $q \in \F[x_1, \ldots, x_k]$. Let $y_1$, \ldots, $y_t \in \F^k$ be a sequence of pairwise distinct points such that $q(y_i) \neq 0$ for all $i \in [t]$. Then, for any $x_1$, \ldots, $x_t \in \F^n$,
		\[
			\dim \set{f \in \Rat_d(k,n;q) \,:\, f(y_1) = x_1, \, \ldots, \, f(y_t) = x_t} \,=\, \dim \Rat_d(k,n;q) - tn \,=\, \left({k+d \choose d} - t\right)n. 
		\]
	\end{lemma}
	\begin{coolproof}
		We can express each $f \in \Rat_d(k,n;q)$ as $f = (p_1/q,\, \ldots, \, p_n/q)$, where $p_1$, \ldots, $p_n$ are polynomials of degree at most $d$. For each $i \in [t]$, write
		$x_i = (x_i(1), \ldots, x_i(n))$ and $y_i = (y_i(1), \ldots, y_i(k))$.
		The conditions $f(y_1) = x_1$, \ldots, $f(y_t) = x_t$ are then equivalent to
		\begin{equation}\label{eq:linear}
			p_j(y_i) = q(y_i) \cdot x_i(j) \qquad \text{for all } i \in [t] \text{ and } j \in [n].
		\end{equation}
		This is a system of $tn$ linear equations in the coefficients of the polynomials $p_1$, \ldots, $p_n$. Therefore, the statement of Lemma~\ref{lemma:injective} is equivalent to the assertion that equations \eqref{eq:linear} are independent. Furthermore, it is enough to establish the independence of the equations corresponding to the same $j \in [n]$, since the equations corresponding to distinct $j$ share no common variables. Since the rows of a matrix $M$ are linearly independent if and only if the linear operator determined by $M$ is surjective, it remains to show that for all $b_1$, \ldots, $b_t \in \F$, there exists a polynomial $p$ of degree at most $d$ such that
		\begin{equation}\label{eq:linear1}
			p(y_i) = b_i \qquad \text{for all } i \in [t].
		\end{equation}
		For each $i \in [t] \setminus \set{1}$, choose some $j_i \in [k]$ such that $y_i(j_i) \neq y_1(j_i)$ and let
		\[
			q_1(v_1, \ldots, v_k) \,\defeq\, \frac{(v_{j_2} - y_2(j_2)) \cdots (v_{j_k} - y_k(j_k))}{(y_1(j_2) - y_2(j_2)) \cdots (y_1(j_k) - y_k(j_k))} \,\in\, \F[v_1, \ldots, v_k].
		\]
		Then $q_1$ is a polynomial of degree $t-1 \leq d$ such that $q_1(y_1) = 1$ and $q_1(y_i) = 0$ for all $i \in [t] \setminus \set{1}$. Similarly, there exist polynomials $q_2$, \ldots, $q_t$ such that $q_i(y_i) = 1$ and $q_j(y_i) = 0$ for $i \neq j$. Clearly, the polynomial $p \defeq b_1 q_1 + \ldots + b_t q_t$ is a solution to \eqref{eq:linear1}, as desired.
	\end{coolproof}

	Recall that we call a set $E \subseteq (\F^n)^t$ injective if for all $(x_1, \ldots, x_t) \in E$, the elements $x_1$, \ldots, $x_t$ are pairwise distinct.
	
	\begin{lemma}\label{lemma:size_of_E}
		Let $t \geq 1$ and let $E \subseteq (\F^n)^t$ be an injective definable set. Suppose that for some $d \geq t-1$ and $q \in \F[x_1, \ldots, x_k] \setminus \set{0}$, we have $\forall^\ast f \in \Rat_d(k,n;q) \, (\dim E[f] \geq k)$. Then $\dim E \geq tn - (t-1)k$.
	\end{lemma}
	\begin{coolproof}
		Consider the set
		\[
			\mathfrak{G} \defeq \set{(f;\, y_1, \ldots, y_t;\, x_1, \ldots, x_t) \in \Rat_d(k,n;q) \times (\F^k)^t \times E \,:\, f(y_1) = x_1, \, \ldots, \, f(y_t) = x_t}.
		\]
		By definition, $(y_1, \ldots, y_t) \in E[f]$ when for some $(x_1, \ldots, x_t) \in E$, we have $(f; y_1, \ldots, y_t; x_1, \ldots, x_t) \in \mathfrak{G}$.
		Since $\forall^\ast f \in \Rat_d(k,n;q) \, (\dim E[f] \geq k)$, this implies that
		\begin{equation}\label{eq:X_lower_bound}
			\dim \mathfrak{G} \,\geq\, \dim \Rat_d(k, n;q) + k.
		\end{equation}
		On the other hand, consider any $(f; y_1, \ldots, y_t; x_1, \ldots, x_t) \in \mathfrak{G}$. By definition, $\set{y_1, \ldots, y_t} \subseteq \dom(f)$, i.e., $q(y_i) \neq 0$ for all $i \in [t]$. Since $E$ is injective, the points $x_1$, \ldots, $x_t$, and hence also $y_1$, \ldots, $y_t$, are pairwise distinct. By Lemma~\ref{lemma:injective}, if we choose any $(x_1, \ldots, x_t) \in E$ and any sequence $y_1$, \ldots, $y_t \in \F^k$ of pairwise distinct points such that $q(y_i) \neq 0$ for all $i \in [t]$, then
		\[
			\dim \set{f \in \Rat_d(k, n;q) \,:\, (f;\, y_1, \ldots, y_t;\, x_1, \ldots, x_t) \in \mathfrak{G}} \,=\, \dim \Rat_d(k, n; q) - tn.
 		\]
 		Therefore,
 		\begin{equation}\label{eq:X_upper_bound}
	 		\dim \mathfrak{G} \,=\, \dim E + tk + \dim \Rat_d(k, n; q) - tn.
 		\end{equation}
 		Combining \eqref{eq:X_lower_bound} and \eqref{eq:X_upper_bound} yields the desired conclusion.
	\end{coolproof}
	
	With Lemma~\ref{lemma:size_of_E} in hand, we derive Proposition~\ref{prop:lower_bound}, restated below for the reader's convenience:
	
	\begin{propcopy}{prop:lower_bound}
		Let $t \geq 1$ and let $E \subseteq (\F^n)^t$ be an injective definable set. If there exist $d \geq t-1$ and a nonzero polynomial $q \in \F[x_1, \ldots, x_k]$ such that, for generic $f \in \Rat_d(k,n;q)$, every definable $E[f]$-independent set $I \subseteq \F^k$ has dimension less than $k$, then $E$ has a $k$-almost dense irreducible component.
	\end{propcopy}
	\begin{coolproof}
		First we show that $E$ itself is $k$-almost dense. Let $\0 \neq S \subseteq [t]$. Since every $(\proj_S E)[f]$-independent set is also $E[f]$-independent, we conclude that for generic $f \in \Rat_d(k,n;q)$, every definable $(\proj_S E)[f]$-independent set $I \subseteq \F^k$ has dimension less than $k$. By Proposition~\ref{prop:dense}, this implies that $\dim (\proj_S E)[f] \geq k$. Then, by Lemma~\ref{lemma:size_of_E}, $\dim \proj_S E \geq |S|n - (|S|-1)k$, as desired.
		
		Now suppose that the components of $E$ are $H_1$, \ldots, $H_s$. Consider any $f \in \Rat_d(k,n;q)$. Notice that $E[f] = H_1[f] \cup \ldots \cup H_s[f]$, and thus if there is no $k$-dimensional definable $E[f]$-independent set, then there is also no $k$-dimensional definable $H_i[f]$-independent set for some $i \in [s]$ \ep{as the intersection of finitely many dense definable sets is dense}. 
		Since the set $\Rat_d(k,n;q)$ is irreducible, there is some $i \in [s]$ such that for generic $f \in \Rat_d(k,n;q)$, every definable $H_i[f]$-independent set has dimension less than $k$. By the above argument, $H_i$ is $k$-almost dense, and we are done.
	\end{coolproof}
	
	Lemma~\ref{lemma:injective} will be used again in \S\S\ref{subsec:induction} and \ref{subsec:finish}.
	
	\subsection{Restrictions}
	
	Fix a nonzero polynomial $q \in \F[x_1, \ldots, x_k]$. Given $g \in \Rat_d(k,n;q)$ and a subset $L \subseteq \L_{<k}(k)$, we use $\Rat_d(k,n;q)[g;L]$ to denote the set of all $f \in \Rat_d(k,n;q)$ such that
	\[
		f \circ \ell = g \circ \ell \qquad \text{for all} \qquad \ell \in L.
	\]
	We view $\Rat_d(k,n;q)[g;L]$ as a definable family of functions, with the evaluation map inherited from $\Rat_d(k,n;q)$, and call it the \emphd{$(g; L)$-restriction} of $\Rat_d(k,n;q)$. Note that, by definition, $\Rat_d(k, n;q)[g;L] \neq \0$ as $g \in \Rat_d(k,n;q)[g;L]$. Since the set $\Rat_d(k,n;q)[g;L]$ is cut out from $\Rat_d(k,n;q)$ by linear equations, it is irreducible.
	
	The next proposition allows us to apply the results of Section \ref{sec:proj} to $\Rat_d(k,n;q)[g;L]$ \ep{it is another reason for the lower bound $d \geq t-1$ in Theorem~\ref{theo:main}}:
	
	\begin{prop}\label{prop:no_cont}
		Fix $q \in \F[x_1, \ldots, x_k]\setminus \set{0}$. Let $g \in \Rat_d(k,n;q)$ and $L \subseteq \L_{<k}(k)$.
		\begin{enumerate}[label=\ep{\itshape{\alph*}}, wide]
			\item\label{item:no_cont_a} If $|L| \leq d$, then the family $\Rat_d(k,n;q)[g;L]$ is comprehensive.
			\item\label{item:no_cont_b} If $|L| < d$, then the family $\Rat_d(k,n;q)[g;L]$ is $(k-1)$-uncontainable. 
		\end{enumerate}
	\end{prop}
	\begin{coolproof}
		\ref{item:no_cont_a} Assume that $|L| \leq d$ and suppose that $\Rat_d(k,n;q)[g;L]$ is not comprehensive. Due to Proposition~\ref{prop:C1toC2}, there is a definable set $C \subset \F^n$ such that $\dim C < n$ but
		\begin{equation}\label{eq:gLcont}
			\forall^\ast f \in \Rat_d(k,n;q)[g;L] \quad \forall^\ast y \in \F^k \qquad f(y) \in C.
		\end{equation}
		We may replace $C$ by its closure and assume that $C$ is closed, i.e., it is the set of common zeros of a family of $n$-variable polynomials. Let $Z \subset \F^k$ be the zero locus of $q$. Since $C$ is closed, the set
		\[
			\set{(f, y) \in \Rat_d(k,n;q)[g;L] \times (\F^k \setminus Z) \,:\, f(y) \in C}
		\]
		is cut out from $\Rat_d(k,n;q)[g;L] \times (\F^k \setminus Z)$ by a family of polynomial equations; in other words, it is relatively closed in $\Rat_d(k,n;q)[g;L] \times (\F^k \setminus Z)$. Since it is dense in $\Rat_d(k,n;q)[g;L] \times (\F^k \setminus Z)$ by \eqref{eq:gLcont}, we can replace generic quantifiers in \eqref{eq:gLcont} by universal ones and conclude that
		\[
			\forall f \in \Rat_d(k,n;q)[g;L] \quad \forall y \in \F^k \setminus Z \qquad f(y) \in C.
		\]
		To arrive at a contradiction, we shall exhibit $f \in \Rat_d(k,n;q)[g;L]$ and $y \in \F^k \setminus Z$ such that $f(y) \not \in C$. Actually, for any given $x \in \F^n$, we will find $f \in \Rat_d(k,n;q)[g;L]$ and $y \in \F^k \setminus Z$ with $f(y) = x$.
		
		Let $L \eqqcolon \set{\ell_1, \ldots, \ell_s}$, where $s \leq d$. Since each $\ell_i$ is defined on a space of dimension less than $k$, we can choose nonconstant affine maps $\phi_i \colon \F^k \to \F$ such that $\im (\ell_i) \subseteq \ker (\phi_i)$. Let \[p \defeq \phi_1 \cdots \phi_s.\] Then $p$ is a nonzero polynomial in $k$ variables of degree $s \leq d$ such that $p \circ \ell_i = 0$ for all $i \in [s]$. For any sequence of coefficients $a_1$, \ldots, $a_n \in \F$, define a rational map $f_{a_1, \ldots, a_n} \in \Rat(k,n)$ via 
		\[
			f_{a_1, \ldots, a_n} \defeq g + \left(\frac{a_1p}{q}, \,\ldots,\, \frac{a_n p}{q}\right).
		\]
		By the choice of $p$, we have $f_{a_1, \ldots, a_n} \in \Rat_d(k,n;q)[g;L]$. Furthermore, if $y \in \F^k \setminus Z$ is such that $p(y) \neq 0$, then, by varying $a_1$, \ldots, $a_n$, we can force $f_{a_1, \ldots, a_n}(y)$ to take an arbitrary value $x \in \F^n$, as desired. 
		
		\ref{item:no_cont_b} Assume that $|L| < d$ and let $\mathfrak{C}$ be a $(k-1)$-container for $\Rat_d(k,n;q)[g;L]$. Then \ref{item:U_capture} becomes
		\[
			\forall^\ast f \in \Rat_d(k,n;q)[g;L] \quad \forall^\ast \ell \in \L(k-1, k) \quad \forall^\ast y \in \F^k \qquad f(y) \in \mathfrak{C}_{f \circ \ell},
 		\]
 		or, equivalently,
 		\begin{equation}\label{eq:U11}
	 		\forall^\ast \ell \in \L(k-1, k) \quad \forall^\ast f \in \Rat_d(k,n;q)[g;L] \quad \forall^\ast y \in \F^k \qquad f(y) \in \mathfrak{C}_{f \circ \ell}.
 		\end{equation}
 		For $\ell \in \L(k-1, k)$, define an equivalence relation $R_\ell$ on $\Rat_d(k,n;q)[g;L]$ by
 		\[
	 		(f, h) \in R_\ell \,\vcentcolon\Longleftrightarrow\, f \circ \ell = h \circ \ell.
 		\]
 		By definition, the $R_\ell$-equivalence class of $f$ is the set $\Rat_d(k,n;q)[f; L \cup \set{\ell}]$. Since $R_\ell$ is cut out from $\Rat_d(k,n;q)[g;L] \times \Rat_d(k,n;q)[g;L]$ by linear equations, it is irreducible, and since $R_\ell$ is reflexive, we have $\proj_1 R_\ell = \proj_2 R_\ell = \Rat_d(k,n;q)[g;L]$. Thus, we may apply implication \ref{item:X} $\Longrightarrow$ \ref{item:E} of Corollary~\ref{corl:equivalence} to rewrite \eqref{eq:U11} as follows:
 		\begin{align}
	 		\forall^\ast \ell \in \L(k-1, k) &\quad \forall^\ast f \in \Rat_d(k,n;q)[g;L] \nonumber\\
	 		&\forall^\ast h \in \Rat_d(k,n;q)[f; L \cup \set{\ell}] \quad \forall^\ast y \in \F^k \qquad h(y) \in \mathfrak{C}_{h \circ \ell}. \label{eq:long}
 		\end{align}
 		Note that if $(f, h) \in R_\ell$, then $h \circ \ell = f \circ \ell$, and hence $\mathfrak{C}_{h \circ \ell} = \mathfrak{C}_{f \circ \ell}$. Also, since $|L| < d$, by part \ref{item:no_cont_a}, the family $\Rat_d(k,n;q)[f; L \cup \set{\ell}]$ is comprehensive. Therefore, we have
 		\begin{align*}
	 		&\forall^\ast h \in \Rat_d(k,n;q)[f; L \cup \set{\ell}] \quad \forall^\ast y \in \F^k \quad h(y) \in \mathfrak{C}_{h \circ \ell} \\
	 		[\mathfrak{C}_{h \circ \ell} = \mathfrak{C}_{f \circ \ell}] \quad \Longleftrightarrow \quad &\forall^\ast h \in \Rat_d(k,n;q)[f; L \cup \set{\ell}] \quad \forall^\ast y \in \F^k \quad h(y) \in \mathfrak{C}_{f \circ \ell} \\
	 		[\text{comprehensiveness}] \quad \Longleftrightarrow \quad &\forall^\ast x \in \F^n \quad x \in \mathfrak{C}_{f \circ \ell}.
 		\end{align*}
 		Hence, \eqref{eq:long} is equivalent to
 		\[
	 		\forall^\ast \ell \in \L(k-1, k) \quad \forall^\ast f \in \Rat_d(k,n;q)[g;L] \quad \forall^\ast x \in \F^n \qquad x \in \mathfrak{C}_{f \circ \ell}, 
 		\]
 		which turns into \ref{item:U_small} after switching the order of the first two quantifiers.
	\end{coolproof}
	
	\subsection{Iterative applications of Theorem~\ref{theo:proj}}\label{subsec:induction}
	
	Let $t \geq 1$ and let $E \subseteq (\F^n)^t$ be a definable set. Let $s \leq t - 1$. Given $f \in \Rat(k,n)$ and $\ell_1$, \ldots, $\ell_s \in \L(r, k)$, define $E(f; \ell_1, \ldots, \ell_s) \subseteq (\F^n)^{t-s}$ as the set of all tuples $(x_{1}, \ldots, x_{t-s}) \in (\F^n)^{t-s}$ such that
	\[
		\exists z_1, \,\ldots,\, z_s \in \F^r \qquad ((f \circ \ell_1)(z_1), \, \ldots, \, (f \circ \ell_s)(z_s), \, x_{1}, \, \ldots, \, x_{t-s}) \,\in\, E.
	\] 
	For $s = 0$, we set $E(f; \0) \defeq E$. It is clear from this definition that for $s \geq 1$,
	\[
		E(f; \ell_1, \ldots, \ell_s) \,=\, E(f; \ell_1, \ldots, \ell_{s-1})(f;\ell_s). 
	\]
	We view $E(f; \ell_1, \ldots, \ell_s)$ as a subset of the product space $(\F^n)^{t-s}$ with the coordinates indexed by $[t-s]$.
	
	\begin{lemma}\label{lemma:induction}
		Let $t \geq 1$ and let $E \subseteq (\F^n)^t$ be an injective $r$-almost dense irreducible definable set. Fix $d \geq t-1$ and a nonzero polynomial $q \in \F[x_1, \ldots, x_k]$. If $k \geq r+1$ and $s \leq t-1$, then, for generic $f \in \Rat_d(k,n;q)$ and $\ell_1$, \ldots, $\ell_s \in \L(r, k)$, the set $E(f; \ell_1, \ldots, \ell_s) \subseteq (\F^n)^{t-s}$ is nonempty and all its irreducible components are $r$-almost dense.
	\end{lemma}
	\begin{coolproof}
		The proof is by induction on $s$. The base case $s=0$ is trivial, so assume $s \geq 1$. We first show that the set $E(f; \ell_1, \ldots, \ell_s)$ itself is $r$-almost dense \ep{and, in particular, nonempty}. The same argument as in the proof of Proposition~\ref{prop:no_cont}\ref{item:no_cont_b} shows that the sequence of quantifiers
		\[
			\forall^\ast f \in \Rat_d(k,n;q) \quad \forall^\ast \ell_1,\, \ldots,\, \ell_s \in \L(r, k)
		\]
		can be replaced by
		\begin{align*}
			\forall^\ast g \in \Rat_d(k,n;q) &\quad \forall^\ast \ell_1,\,\ldots,\, \ell_{s-1} \in \L(r, k) \\
			&\forall^\ast f \in \Rat_d(k,n;q)[g; \set{\ell_1, \ldots, \ell_{s-1}}] \quad \forall^\ast \ell_s \in \L(r,k).
		\end{align*}
		Consider generic $g \in \Rat_d(k,n;q)$ and $\ell_1$, \ldots, $\ell_{s-1} \in \L(r,k)$ and let
		\[
			\calF \defeq \Rat_d(k,n;q)[g; \set{\ell_1, \ldots, \ell_{s-1}}] \qquad \text{and} \qquad E' \defeq E(g; \ell_1, \ldots, \ell_{s-1}).
		\]
		Notice that if $f \in \calF$, then $E(f; \ell_1, \ldots, \ell_{s-1}) = E'$, and thus for any $\ell_s \in \L(r,k)$,
		\[
			E(f; \ell_1, \ldots, \ell_s) = E'(f;\ell_s).
		\]
		By the inductive assumption, $E' \neq \0$ and every irreducible component of $E'$ is $r$-almost dense. Let $H$ be any component of $E'$ \ep{so $H$ is $r$-almost dense} and let $\0 \neq S \subseteq [t-s]$. Set $S' \defeq \set{i + 1 \,:\, i \in S} \subset [t-s+1]$. Then for each $\ell_s \in \L(r,k)$, we have
		\[
			\proj_S (E(f; \ell_1, \ldots, \ell_s)) \,=\, \proj_S (E'(f; \ell_s)) \,\supseteq\, \proj_S (H(f; \ell_s)) \,=\, (\proj_{\set{1} \cup S'} H)(f; \ell_s).
		\]
		Define $A \defeq \proj_{\set{1} \cup S'} H$. We can view $A$ as a subset of the product space $\F^n \times (\F^n)^{S'}$. Then
		\[
			A(f; \ell_s) \,=\, \proj_2 (A_{f \circ \ell_s}),
		\] 
		where $A_{f \circ \ell_s}$ denotes the fiber of $A$ over $f \circ \ell_s$ \ep{as defined in \S\ref{subsec:proj}}.
		Since $H$ is $r$-almost dense, we have \[\dim \proj_1 A \,=\, \dim \proj_1 H \,=\, n,\] i.e., the set $A$ is dominant. Since $r < k$ and $s-1 \leq t-2 < d$, Proposition~\ref{prop:no_cont}\ref{item:no_cont_b} implies that the family $\calF$ is $(k-1)$-uncontainable, and, by Proposition~\ref{prop:reduce}, the family $\calF \otimes \L(r,k)$ is $(r-1)$-uncontainable. Therefore, since $A$ is irreducible, we can apply Lemma~\ref{lemma:dominant} and Theorem~\ref{theo:proj} to conclude that, for generic $f \in \calF$ and $\ell_s \in \L(r,k)$,
		\begin{equation}\label{eq:appl_of_theo}
			\dim\proj_S (E(f; \ell_1, \ldots, \ell_s)) \,\geq\, \dim \proj_2 (A_{f \circ \ell_s}) \,\geq\, \min \set{\dim A - n + r,\, \dim \proj_2 A}.
		\end{equation}
		Observe that, since $H$ is $r$-almost dense,
		\[
			\dim A \,=\, \dim \proj_{\set{1} \cup S'} H \,\geq\, |\set{1} \cup S'| n - (|\set{1} \cup S'| - 1)r \,=\, (|S|+1)n - |S|r, 
		\]
		so $\dim A -n + r \geq |S|n - (|S|-1)r$.
		Similarly,
		\[
			\dim \proj_2 A \,=\, \dim \proj_{S'} H \,\geq\, |S'|n - (|S'|-1)r \,=\, |S|n - (|S|-1)r. 
		\]
		Thus, \eqref{eq:appl_of_theo} yields \[\dim\proj_S (E(f; \ell_1, \ldots, \ell_s)) \,\geq\, |S|n - (|S|-1)r,\] and hence $E(f; \ell_1, \ldots, \ell_s)$ is $r$-almost dense, as claimed.
		
		To prove that every component of $E(f; \ell_1, \ldots, \ell_s)$ is $r$-almost dense, let $\mathfrak{X}$ be the set of all tuples
		\[
			(f;\,\ell_1, \ldots, \ell_s;\, x_1, \ldots, x_{t-s}) \,\in\, \Rat_d(k,n;q) \times (\L(r,k))^s \times (\F^n)^{t-s}
		\]
		such that $(x_1, \ldots, x_{t-s}) \in E(f; \ell_1, \ldots, \ell_s)$.
		
		\stepcounter{ForClaims} \renewcommand{\theForClaims}{\ref{lemma:induction}}
		\begin{claim}\label{claim:big_irred_induct}
			The set $\mathfrak{X}$ is irreducible.
		\end{claim}
		\begin{claimproof}
			Writing $f = (p_1/q, \ldots, p_n/q)$, let $\mathfrak{Z}$ be the set of all tuples
			\[
				(f;\, \ell_1, \ldots, \ell_s;\, z_1, \ldots, z_s; \, y_1, \ldots, y_s;\, x_1, \ldots, x_t) \,\in\, \Rat_d(k,n;q) \times (\L(r,k))^s \times (\F^r)^s \times (\F^k)^s \times (\F^n)^t
			\]
			such that $\ell_i(z_i) = y_i$ and $(p_1(y_i),\ldots,p_n(y_i)) = q(y_i) \cdot x_i$ for all $i \in [s]$. Let $D$ denote the set of all tuples $(y_1, \ldots, y_s) \in (\F^k)^s$ of pairwise distinct points such that $q(y_i) \neq 0$ for all $i \in [s]$, and let
			\[
				\mathfrak{Y} \,\defeq\, \mathfrak{Z} \cap (\Rat_d(k,n;q) \times (\L(r,k))^s \times (\F^r)^s \times D \times E).
			\]
			Note that $D$ is irreducible, as $\cl{D} = (\F^k)^s$. Since $E$ is injective, $\mathfrak{X}$ is the image of $\mathfrak{Y}$ under the projection
			\[
				(f;\, \ell_1, \ldots, \ell_s;\, z_1, \ldots, z_s; \, y_1, \ldots, y_s;\, x_1, \ldots, x_t) \,\mapsto\, (f;\, \ell_1, \ldots, \ell_s;\, x_{s+1}, \ldots, x_t),
			\]
			so it suffices to prove that $\mathfrak{Y}$ is irreducible. The defining equations for $\mathfrak{Z}$ are linear in the coefficients of $p_1$, \ldots, $p_n$ and $\ell_1$, \ldots, $\ell_s$, so we may use Corollary~\ref{corl:linear_irred}. The set $(\F^r)^s \times D \times E$ is irreducible, so we just need to show that for all $(z_1, \ldots, z_s; y_1, \ldots, y_s; x_1, \ldots, x_t) \in (\F^r)^s \times D \times E$, the dimension of the set
			\[
				\set{(f;\, \ell_1, \ldots, \ell_s) \,\in\, \Rat_d(k,n;q) \times (\L(r,k))^s \,:\, (f;\, \ell_1, \ldots, \ell_s;\, z_1, \ldots, z_s; \, y_1, \ldots, y_s;\, x_1, \ldots, x_t) \,\in\, \mathfrak{Z}}
			\]
			is the same. But since $(y_1, \ldots, y_s) \in D$, this is indeed the case by Lemma~\ref{lemma:injective}, according to which the dimension of this set is equal to $\dim \Rat_d(k,n;q) - sn + s(\dim \L(r,k) - k)$.
		\end{claimproof}
		
		We can now finish the proof of Lemma~\ref{lemma:induction}. Take any $\0 \neq S \subseteq [t - s]$ and let $\mathfrak{X}_S \subseteq \mathfrak{X}$ be the set such that for all $f \in \Rat_d(k, n; q)$ and $\ell_1$, \ldots, $\ell_s \in \L(r, k)$, the fiber of $\mathfrak{X}_S$ over $(f; \ell_1, \ldots, \ell_s)$ is the union of all the irreducible components $H$ of $E(f; \ell_1, \ldots, \ell_s)$ with
		\begin{equation}\label{eq:H}
			\dim \proj_S H \geq |S|n - (|S|-1)r.
		\end{equation}
		It follows from the results of \S\ref{subsec:model_theory} that $\mathfrak{X}_S$ is definable. We already know that for generic $f$, $\ell_1$, \ldots, $\ell_s$, the set $E(f; \ell_1, \ldots, \ell_s)$ is $r$-almost dense, and hence $\dim \proj_S (E(f;\ell_1, \ldots, \ell_s)) \geq |S|n - (|S|-1)r$, which means that $E(f;\ell_1, \ldots, \ell_s)$ has a component $H$ satisfying \eqref{eq:H}, i.e., the fiber of $\mathfrak{X}_S$ over $(f;\ell_1, \ldots, \ell_s)$ is nonempty. Since $\mathfrak{X}$ is irreducible, Corollary~\ref{corl:split} implies that for generic $f$, $\ell_1$, \ldots, $\ell_s$, the fiber of $\mathfrak{X}_S$ over $(f;\ell_1, \ldots, \ell_s)$ must be equal to $E(f;\ell_1, \ldots, \ell_s)$. In other words, every irreducible component $H$ of $E(f;\ell_1, \ldots, \ell_s)$ satisfies \eqref{eq:H}, as desired.
	\end{coolproof}
	
	Applying Lemma~\ref{lemma:induction} with $s = t-1$ yields the following:
	
	\begin{corl}\label{corl:prints}
		Let $t \geq 1$ and let $E \subseteq (\F^n)^t$ be an injective $r$-almost dense irreducible definable set. Fix $d \geq t-1$ and a nonzero polynomial $q \in \F[x_1, \ldots, x_k]$. If $k \geq r+1$, then
		\begin{align*}
		\forall^\ast f \in \Rat_d(k, n; q)&\quad \forall^\ast \ell_1, \,\ldots,\, \ell_{t-1} \in \L(r, k)\quad \forall^\ast y \in \F^k \\
		&\exists z_1, \,\ldots,\, z_{t-1} \in \F^r \qquad (\ell_1(z_1),\, \ldots,\, \ell_{t-1}(z_{t-1}),\, y) \in E[f].
		\end{align*}
	\end{corl}
	\begin{coolproof}
		We again observe that the sequence of quantifiers
		\[
			\forall^\ast f \in \Rat_d(k, n; q)\quad \forall^\ast \ell_1, \,\ldots,\, \ell_{t-1} \in \L(r, k)
		\]
		can be replaced by
		\[
			\forall^\ast g \in \Rat_d(k, n; q)\quad \forall^\ast \ell_1, \,\ldots,\, \ell_{t-1} \in \L(r, k) \quad \forall^\ast f \in \Rat_d(k,n;q)[g; \set{\ell_1, \ldots, \ell_{t-1}}].
		\]
		Consider generic $g \in \Rat_d(k,n;q)$ and $\ell_1$, \ldots, $\ell_{t-1} \in \L(r,k)$ and let
		\[
			\calF \defeq \Rat_d(k,n;q)[g; \set{\ell_1, \ldots, \ell_{t-1}}] \qquad \text{and} \qquad C \defeq E(g; \ell_1, \ldots, \ell_{t-1}) \subseteq \F^n.
		\]
		Note that for all $f \in \calF$, $E(f;\ell_1, \ldots, \ell_{t-1}) = C$. Since $d \geq t-1$, Proposition~\ref{prop:no_cont}\ref{item:no_cont_a} shows that $\calF$ is comprehensive, and hence we obtain
		\begin{align*}
			&\forall^\ast f \in \calF \quad \forall^\ast y \in \F^k \quad \exists z_1, \,\ldots,\, z_{t-1} \in \F^r \quad (\ell_1(z_1),\, \ldots,\, \ell_{t-1}(z_{t-1}),\, y) \in E[f] \\
			\Longleftrightarrow \quad & \forall^\ast f \in \calF \quad \forall^\ast y \in \F^k \quad f(y) \in E(f; \ell_1, \ldots, \ell_{t-1}) \\
			\Longleftrightarrow \quad & \forall^\ast f \in \calF \quad \forall^\ast y \in \F^k \quad f(y) \in C \\
			\Longleftrightarrow \quad & \forall^\ast x \in \F^n \quad x \in C \\
			\Longleftrightarrow \quad & \dim C = n.
		\end{align*}
		But $\dim C = n$ is precisely the conclusion of Lemma~\ref{lemma:induction} for $s = t-1$, so we are done.
	\end{coolproof}
	
	\subsection{Finishing the proof}\label{subsec:finish}
	
	We are finally ready to prove Theorem~\ref{theo:main}. For the reader's convenience, we state it again below:
	
	\begin{theocopy}{theo:main}
		Let $t \geq 1$ and let $E \subseteq (\F^n)^t$ be an $r$-almost dense irreducible definable set. Fix $d \geq t-1$ and a nonzero polynomial $q \in \F[x_1, \ldots, x_k]$. If $k \geq r + 1$, then, for generic $f \in \Rat_d(k,n;q)$, the following holds:
		
		\smallskip
		
		Every definable $E[f]$-independent set $I \subseteq \F^k$ has dimension less than $k$. Furthermore, if $E$ is injective, then every irreducible component of $E[f]$ is $r$-almost dense.
	\end{theocopy}
	\begin{coolproof}
		Let us first assume that $E$ is injective. As in the proof of Lemma~\ref{lemma:induction}, we start by showing that for generic $f \in \Rat_d(k,n;q)$, the set $E[f]$ itself is $r$-almost dense \ep{hence nonempty}. Consider any $\0 \neq S \subseteq [t]$. For concreteness, we may assume that $S = [s] \cup \set{t}$ for some $s \leq t-1$. It follows from Corollary~\ref{corl:prints} that a generic map $f \in \Rat_d(k,n;q)$ satisfies
		\begin{align*}
			&\quad \forall^\ast \ell_1, \,\ldots,\, \ell_s \in \L(r, k)\quad \forall^\ast y \in \F^k \\
		&\exists z_1, \,\ldots,\, z_s \in \F^r \quad \exists y_{s+1}, \,\ldots,\, y_{t-1} \in \F^k \qquad (\ell_1(z_1),\, \ldots,\, \ell_{s}(z_{s}),\, y_{s+1}, \, \ldots, \, y_{t-1},\, y) \in E[f],
		\end{align*}
		which can be rewritten as
		\begin{align}
		&\quad \forall^\ast \ell_1, \,\ldots,\, \ell_s \in \L(r, k)\quad \forall^\ast y \in \F^k \nonumber\\
		&\exists z_1, \,\ldots,\, z_s \in \F^r \qquad (\ell_1(z_1),\, \ldots,\, \ell_{s}(z_{s}),\, y) \in \proj_S (E[f]). \label{eq:counting}
		\end{align}
		Let $\mathfrak{G}$ be the set of all tuples \[(\ell_1, \ldots, \ell_s;\,  y_1, \ldots, y_{s+1};\,  z_1, \ldots, z_s) \,\in\, (\L(r,k))^s \times \proj_S(E[f]) \times (\F^r)^s\] such that $\ell_i(z_i) = y_i$ for all $i \in [s]$.
		From \eqref{eq:counting}, we conclude that
		\begin{equation}\label{eq:main_X_lower}
			\dim \mathfrak{G} \,\geq\, s \dim \L(r,k) + k.
		\end{equation}
		On the other hand, for all $y \in \F^k$ and $z \in \F^r$, we have $\dim \set{\ell \in \L(r,k) \,:\, \ell(z) = y} = \dim \L(r,k) - k$ by Lemma~\ref{lemma:injective}. Hence, for all $(y_1, \ldots, y_s, y_{s+1}) \in \proj_S(E[f])$ and $z_1$, \ldots, $z_s \in \F^r$, the set of all tuples $(\ell_1, \ldots, \ell_s) \in (\L(r,k))^s$ with $\ell_1(z_1) = y_1$, \ldots, $\ell_s(z_s) = y_s$ has dimension $s(\dim\L(r,k) - k)$, so
		\begin{equation}\label{eq:main_X_upper}
			\dim \mathfrak{G} \,=\, \dim \proj_S (E[f]) + sr + s(\dim\L(r,k) - k).
		\end{equation}
		Comparing \eqref{eq:main_X_lower} and \eqref{eq:main_X_upper}, we obtain
		\[
			\dim \proj_S (E[f]) \,\geq\, (s+1)k - sr \,=\, |S|k - (|S|-1)r,
		\]
		as desired.
		
		To deduce that every irreducible component of $E[f]$ is $r$-almost dense, we use the same trick as in the proof of Lemma~\ref{lemma:induction}. Define
		\[
			\mathfrak{X} \defeq \set{(f;\, y_1, \ldots, y_t) \in \Rat_d(k,n;q) \times (\F^k)^t \,:\, (y_1, \ldots, y_t) \in E[f]}.
		\]
		
		\stepcounter{ForClaims} \renewcommand{\theForClaims}{\ref{theo:main}}
		\begin{claim}\label{claim:irred_again}
			The set $\mathfrak{X}$ is irreducible.
		\end{claim}
		\begin{claimproof}
			The argument is analogous to the proof of Claim~\ref{claim:big_irred_induct}. Writing $f = (p_1/q, \ldots, p_n/q)$, let $\mathfrak{Z}$ be the set of all tuples
			\[
			(f;\, y_1, \ldots, y_t;\, x_1, \ldots, x_t) \,\in\, \Rat_d(k,n;q) \times (\F^k)^t \times (\F^n)^t
			\]
			such that $(p_1(y_i),\ldots,p_n(y_i)) = q(y_i) \cdot x_i$ for all $i \in [t]$. Let $D$ denote the set of all tuples $(y_1, \ldots, y_t) \in (\F^k)^t$ of pairwise distinct points such that $q(y_i) \neq 0$ for all $i \in [t]$, and let
			\[
			\mathfrak{Y} \,\defeq\, \mathfrak{Z} \cap (\Rat_d(k,n;q) \times D \times E).
			\]
			Since $E$ is injective, $\mathfrak{X}$ is the image of $\mathfrak{Y}$ under the projection
			\[
			(f; \, y_1, \ldots, y_t;\, x_1, \ldots, x_t) \,\mapsto\, (f;\, y_1, \ldots, y_t),
			\]
			so it suffices to prove that $\mathfrak{Y}$ is irreducible. The defining equations for $\mathfrak{Z}$ are linear in the coefficients of $p_1$, \ldots, $p_n$, so we may use Corollary~\ref{corl:linear_irred}. The set $D \times E$ is irreducible, so we just need to show that for all $(y_1, \ldots, y_t; x_1, \ldots, x_t) \in D \times E$, the dimension of the set
			\[
			\set{f \in \Rat_d(k,n;q) \,:\, (f;\, y_1, \ldots, y_t;\, x_1, \ldots, x_t) \,\in\, \mathfrak{Z}}
			\]
			is the same. But since $(y_1, \ldots, y_t) \in D$ and $d \geq t-1$, this is the case by Lemma~\ref{lemma:injective}, according to which the dimension of this set is equal to $\dim \Rat_d(k,n;q) - tn$.
		\end{claimproof}
		
		Treating each $f \in \Rat_d(k,n;q)$ as a point in the space $\Rat_d(k,n;q)$, we may consider the fiber $\mathfrak{X}_f$ of $\mathfrak{X}$ over $f$. By definition,
		\[
			\mathfrak{X}_f \,=\, \set{(y_1, \ldots, y_t) \in (\F^k)^t \,:\, (f; y_1, \ldots, y_t) \in \mathfrak{X}} \,=\, E[f].
		\]
		Now take any $\0 \neq S \subseteq [t]$ and let $\mathfrak{X}_S \subseteq \mathfrak{X}$ be the set such that for all $f \in \Rat_d(k, n; q)$, the fiber \[
			(\mathfrak{X}_S)_f \,=\, \set{(y_1, \ldots, y_t) \in (\F^k)^t \,:\, (f; y_1, \ldots, y_t) \in \mathfrak{X}_S} \,\subseteq\, E[f]
		\]
		is the union of all the irreducible components $H$ of $E[f]$ with $\dim \proj_S H \geq |S|k - (|S|-1)r$. Then $\mathfrak{X}_S$ is definable, and we have shown that for generic $f \in \Rat_d(k,n;q)$, $(\mathfrak{X}_S)_f \neq \0$. Since $\mathfrak{X}$ is irreducible, Corollary~\ref{corl:split} implies that for generic $f \in \Rat_d(k,n;q)$, $(\mathfrak{X}_S)_f = E[f]$, i.e., every component $H$ of $E[f]$ satisfies $\dim \proj_S H \geq |S|k - (|S|-1)r$, as desired. Recall that since, for generic $f \in \Rat_d(k,n;q)$, the set $E[f]$ is nonempty and all its components are $r$-almost dense, Observation~\ref{obs:ind} yields that there is no $k$-dimensional definable $E[f]$-independent set.
		
		Now suppose that $E$ is not necessarily injective. We need to show that for generic $f \in \Rat_d(k,n;q)$, every definable $E[f]$-independent set has dimension less than $k$. Towards a contradiction, let $E$ be a counterexample with the smallest value of $t$. Let
		\[
			E_0 \defeq \set{(x_1, \ldots, x_t) \in E \,:\, x_1, \, \ldots, \, x_t \text{ are pairwise distinct}},
		\]
		and for $1 \leq i < j \leq t$, define
		\[
			E_{ij} \defeq \set{(x_1, \ldots, x_t) \in E \,:\, x_i = x_j}.
		\]
		Since $E$ is irreducible, at least one of the sets $E_0$, $E_{ij}$, $1 \leq i < j \leq t$, is dense in $E$. If $\cl{E_0} = \cl{E}$, then the set $E_0$ is irreducible and, by Corollary~\ref{corl:proj_cl}, $r$-almost dense. Since $E_0$ is injective by definition, we conclude that for generic $f \in \Rat_d(k,n;q)$, the set $E_0[f]$ does not admit a $k$-dimensional definable independent set, and hence the same is true for $E[f] \supseteq E_0[f]$. If, on the other hand, $\cl{E_{ij}} = \cl{E}$ for some $1 \leq i < j \leq t$, then the set $E_{ij}$ is irreducible and $r$-almost dense. Let $S \defeq [t] \setminus \set{j}$ and $E' \defeq \proj_S E_{ij}$. Then $E'$ is also irreducible and $r$-almost dense, so, by the minimality of $t$, for generic $f \in \Rat_d(k,n;q)$, there is no $k$-dimensional definable $E'[f]$-independent set. But a set $I \subseteq \F^k$ is $E'[f]$-independent if and only if it is $E_{ij}[f]$-independent. Hence, the set $E_{ij}[f]$, and thus also $E[f] \supseteq E_{ij}[f]$, does not admit a $k$-dimensional definable independent set, and the proof is complete.
	\end{coolproof}

	\section{Further directions}\label{sec:further}
	
	In this paper we worked with hypergraphs definable in an algebraically closed field $\F$. A natural next step would be to study hypergraphs definable in strongly minimal structures, with Morley rank assuming the role of dimension. 
	It seems especially promising to look at strongly minimal structures with the so-called \emph{definable multiplicity property}, which was isolated by Hrushovski in \cite{Hru}, as they satisfy natural analogs of many basic facts stated in Section~\ref{sec:prelim}. While it appears likely that most our arguments could be extended to this more general setting without too much difficulty, at certain places we invoke properties that are very special to algebraically closed fields. In particular, the proofs of Claims~\ref{claim:big_irred_induct} and \ref{claim:irred_again} use Corollary \ref{corl:linear_irred} and Lemma \ref{lemma:same_fibers}, which ultimately rely on $\mathbb{P}^n$ being a complete variety---a fact that has no obvious analog in arbitrary strongly minimal structures. Nevertheless, we suspect that at least the results of Section~\ref{sec:proj} should have their counterparts in the strongly minimal setting.
	
	Another structure of interest is the real field $\R$. There our arguments cannot be easily adapted, because, in contrast to algebraically closed fields, $\R$ lacks a well-behaved notion of irreducibility. In particular, the ``size'' of an $\R$-definable set is characterized not only by its dimension, but also by its \emph{measure}. It makes sense, therefore, to define the \emphd{independence ratio} $\alpha(E)$ of an $\R$-definable hypergraph $E \subseteq ([0,1]^n)^t$ as $\alpha(E) \defeq \sup_I \lambda_n(I)$, where $\lambda_n$ is the $n$-dimensional Lebesgue measure and the supremum is taken over all $\R$-definable $E$-independent sets $I \subseteq [0,1]^n$. \ep{We are considering hypergraphs on $[0,1]^n$ to make the total measure of the vertex set $1$.}
	The goal is to isolate the properties of $E$ which guarantee that, for a generic $\R$-definable map $f \colon [0,1]^k \to [0,1]^n$, the induced subhypergraph $E[f] \subseteq ([0,1]^k)^t$ satisfies $\alpha(E[f]) \leq \alpha(E)$. A similar problem can also be studied over finite fields $\F_q$, using the \ep{normalized} counting measure instead of the Lebesgue measure.
	
	\subsubsection*{Acknowledgments}
	
	We are very grateful to the anonymous referee for carefully reading the manuscript and providing many useful comments and suggestions.
	
	{\renewcommand{\markboth}[2]{}
		\printbibliography}

\end{document}